\documentclass[conference,onecolumn]{IEEEtran}
\IEEEoverridecommandlockouts
\usepackage{amsmath,amssymb,amsfonts}
\usepackage{graphicx}
	
\usepackage{booktabs}
\usepackage{multirow}

\usepackage{subfigure}
\usepackage{amssymb,amsmath,amsfonts}
\usepackage{multirow}

\usepackage{algorithm}
\usepackage[noend]{algpseudocode}

\DeclareMathAlphabet{\mathcal}{OMS}{cmsy}{m}{n}

\algdef{SE}[DOWHILE]{Do}{doWhile}{\algorithmicdo}[1]{\algorithmicwhile\ #1}%

\newcommand*{\Scale}[2][4]{\scalebox{#1}{$#2$}}%

\begin{document}

\title{{\huge{Adaptive Model Refinement with Batch Bayesian Sampling for Optimization of Bio-inspired Flow Tailoring}}
}

\author{\IEEEauthorblockN{Payam Ghassemi\IEEEauthorrefmark{1}, Sumeet Sanjay Lulekar\IEEEauthorrefmark{2}, and
Souma Chowdhury\IEEEauthorrefmark{3}}
\IEEEauthorblockA{\textit{Department of Mechanical and Aerospace Engineering} \\
\textit{University at Buffalo}\\
Buffalo, NY, 14260\\
Email: \IEEEauthorrefmark{1}payamgha@buffalo.edu,\IEEEauthorrefmark{2}sumeetsa@buffalo.edu,
\IEEEauthorrefmark{3}soumacho@buffalo.edu}}

\maketitle
\begin{abstract}\label{sec:Abst}
This paper presents an advancement to an approach for model-independent surrogate-based optimization with adaptive batch sampling, known as Adaptive Model Refinement (AMR). While the original AMR method provides unique decisions with regards to ``when" to sample and ``how many" samples to add (to preserve the credibility of the optimization search process), it did not provide specific direction towards ``where" to sample in the design variable space. This paper thus introduces the capability to identify optimum location to add new samples. The location of the infill points is decided by integrating a Gaussian Process-based criteria (``q-EI"), adopted from Bayesian optimization. The consideration of a penalization term to mitigate interaction among samples (in a batch) is crucial to effective integration of the q-EI criteria into AMR. The new AMR method, called AMR with Penalized Batch Bayesian Sampling (AMR-PBS) is tested on benchmark functions, demonstrating better performance compared to Bayesian EGO. In addition, it is successfully applied to design surface riblets for bio-inspired passive flow control (where high-fidelity samples are given by costly RANS CFD simulations), leading to a 10\% drag reduction over the corresponding baseline (i.e., riblet-free aerodynamic surface).
\end{abstract}

\begin{IEEEkeywords}
Variable-Fidelity Optimization, Surrogate-Based Optimization, Bayesian Optimization, Bio-inspired Fluid Dynamics, Passive Flow Control, Riblets
\end{IEEEkeywords}

\section{Introduction} \label{sec:intro}
Optimizing complex systems often involves computationally expensive simulations (e.g., CFD) to evaluate system behavior and estimate quantities of interest. While surrogate models and surrogate based optimization \cite{simpson2008design,jin2011surrogate,fernandez2016review} provides a tractable alternative, in their native form they tend to compromise on the fidelity of the optimization process. Specifically, low-fidelity surrogate models often mislead the search process during optimization, leading to sub-optimal or even infeasible solutions. \emph{Variable-fidelity optimization} approaches seek to address these issues and offer attractive trade-offs between computational efficiency and fidelity of the optimal solutions obtained, that is, the ability to quickly arrive at optimal solutions that can be relied upon. In doing so, these methods must answer the following questions depending on the interplay of model uncertainty and function improvement over iterations -- \textit{when} (during optimization) and \textit{where} (in the design space) to re-sample in order to refine the model on the fly? There exist very few frameworks that allow flexible choice of a variety of surrogate models (i.e., model independent) and optimization processes, while also answering both these questions. In this paper, we present important fundamental extension of one such method, known as Adaptive Model Refinement \cite{mehmani2015variable,mehmani2015adaptive}; the new advancements extend the capability of the original method by providing a mechanism to decide the optimal location of batches of new samples, allowing effective utilization of resources. This multifidelity optimization framework is first tested on benchmark functions and then applied to design ribleted 3D airfoil surfaces for bio-inspired flow tailoring. The remaining portion of this introduction section provides a brief review of variable fidelity optimization in general, batch sampling methods in specific, and converges on the objectives of this paper. 

\subsection{Variable Fidelity Optimization} \label{ssec:VFO}
Major surrogate modeling methods that have been used in surrogate based optimization include polynomial response surfaces \cite{RJin-02}, Kriging \cite{Simpson_01,Forrester_Response_review_09}, moving least square \cite{Choi-01,Toropov-05}, radial basis functions (RBF) \cite{Hardy_RBF_71}, support vector regression (SVR) \cite{Clarke05}, artificial neural networks \cite{Yegnanarayana_ANN_04} and hybrid surrogate models~\cite{Zhang_SMO_11}.
In \emph{Variable fidelity optimization} approaches, model management strategies, e.g., model selection, switching and/or refinement, adaptively integrate models of different levels of fidelity and (computational) cost into the optimization process. 
Different model management approaches have been reported in the literature for integrating low-fidelity models within optimization algorithms. One class of model management strategies is developed based on the \emph{Trust-Region} concept~\cite{Booker_1999_VM_TR,Alexandrov_00,Rodriguez_01,marduel2006variableTrustRegion,robinson2008TrustRegionsurrogate}.
The \emph{Trust-Region} methods seek the agreement of the function and its gradient values in the low-fidelity model with those estimated in the high-fidelity model. However, these techniques may not be directly applicable in problems where gradients are expensive to evaluate, or where zero-order algorithms are being used for optimization. 

Another class of insitu model management strategies developed for non-physics-based low-fidelity models (aka. surrogate models) involves adding infill points where additional evaluations of the high-fidelity model or experiment are desired to be performed.
In these approaches, infill points are generally added in (i) the region where the optimum is located (local exploitation), and/or (ii) the entire design space in order to improve the global accuracy of the surrogate (global exploration)\cite{Keane_05,Forrester_book_08,Sugiyama-06}. Infill points can be added in a fully sequential (one-at-a-time), or a batch sequential manner. Various criteria exist for determining the locations of the infill points including (i) index-based criteria (e.g., (integrated- and maximum) Mean Square Error (MSE) and maximum entropy criteria) and (ii) distance-based criteria (e.g., Euclidean distance, Mahalonobis distance, and Weighted distance criteria) ~\cite{Jones-98,Williams-11, Kean-01Sta,Booker-99,Audet-00,Rai-06,SBO-Romero-2011}. Jones et al. in 1998~\cite{Jones-98} developed a well-known model management strategy that is based on an \emph{expected improvement (EI)} criterion, and is called \emph{efficient global optimization (EGO)}, essentially a type of Bayesian Optimization approach. Many variations of this method exist in the literature~\cite{hu2015mixed,viana2013efficient,jeong2005efficient,tajbakhsh2013fully}, detailed review of which is outside the scope of this paper. The EGO method and its variations, however, are generally applicable only to surrogate based optimization based on Gaussian process models. 

Jin et al.~\cite{Jin-02} reviewed different one-at-a-time sequential sampling criteria, and illustrated their clear benefits over single-stage methods~\cite{McKay_79,Sacks-89}. Loeppky et al. \cite{Loeppky-10} explored different batch and one-at-a-time sequential criteria for the Gaussian process model, and pointed out potential cost advantages of batch sampling approaches over one-at-a-time augmentation. One major practical advantage of batch sampling is the ability to better exploit parallel computing resources. A more fundamental advantage of batch sampling is tied with the flexible refinement needs to ensure a reliable optimization progress, which may not necessarily be met with a single sample addition. However batch sampling also brings with it new challenges with regards to computational tractability of modeling the potential interactions among the samples in the batch. A brief review of batch sampling methods specific to Bayesian Optimization, converging on the concept that we build upon, is given next. 

\subsection{Batch Sampling in Optimization} \label{ssec:BATCH}
Given the growing popularity of Bayesian optimization (BO) in solving high-dimensional optimization problems and for hyper-parameter tuning of machine learning architectures, and benefits of being able to add samples in a batch as opposed to the one-at-a-time approach in classical BO, a number of batch Bayesian optimization approaches have emerged in recent years. In the work by Ginsbourger et al. \cite{ginsbourger2008multi}, a multi-point expected improvement (qEI) is generalized to the batch setting, by using the qEI acquisition function for batches of $q$ points. However, it is challenging to identify the points that jointly maximize qEI since the computational cost of evaluating the function and its derivative scales poorly with increasing $q$. Later efforts have suggested heuristics for approximating the qEI computation \cite{snoek2012practical, chevalier2013fast, wang2016parallel}. One such method is the Constant Liar approach ~\cite{ginsbourger2010kriging}, which is a sequential batch building method, that iteratively adds a point that maximizes the single point acquisition function. This approach assumes that evaluating this point will return a particular constant `lie' value, which will temporarily augment the model training set with synthetic values and refitting the GP.

In a recent work, González et al.~\cite{gonzalez2016batch} proposed an alternative batching method by approximating the repulsive effect when batching. Under a Gaussian Process prior, target values of nearby points in sample space are expected to be highly correlated. Thus, when choosing a batch of samples, we may wish for the batch members to be sufficiently far apart to maximize the information gained. To do this, the authors propose the Local Penalization method that sequentially assembles batches of samples by successively penalizing the acquisition function around points previously selected, using a penalization radius based on the estimated Lipschitz constant of the acquisition function surface. \textit{We extend this method for usage in determining optimum locations (i.e., with qEI) of the new batch of samples during refinement events invoked by the AMR method \cite{mehmani2015variable}, where the batch size is given by AMR request to meet the needs of balancing model uncertainty with fitness improvement rate.} In doing so, we integrate a lower-level BO search instance (for the purpose of sampling) within the top level AMR-based optimization, where the latter typically uses a population based algorithm such as Particle Swarm Optimization \cite{Souma2013MDPSO}. This integration is done such that the main advantages of AMR, w.r.t. surrogate model-independence, exploitation of powerful (population-based) global optimizers, and ability to decide when to refine (important for computational efficiency), are retained.

\subsection{Objectives of this Paper}
The primary objectives and associated contributions of this paper are given below: 
\begin{enumerate}
\item Integrate, into the Adaptive Model Refinement method, a request-adaptive batch sampling approach that maximizes the EI while mitigating interactions among samples (within the batch) in a tractable manner. The central contribution lies in the ability (thus enabled) to maximally exploit infill sampling resources while preserving the original benefits (e.g., model independence) of the AMR framework. 
\item Test the upgraded AMR method, called \textit{AMR with Penalized Batch-Bayesian Sampling} or AMR-PBS, on different benchmark functions and compare its performance with that of the Bayesian Efficient Global Optimization (BEGO) method.
\item Apply the new AMR-PBS method to perform multi-fidelity optimization of bio-inspired surface riblets on 3D airfoil sections, to achieve maximum drag reduction. This objective lends itself to the practical contribution of this paper -- namely provisions evidence of the effectiveness of the upgraded AMR method is solving complex multi-fidelity optimization problems.

 \end{enumerate}
The remainder of the paper is structured as follows. The next section provides an overview of the Adaptive Model Refinement method, followed by a detailed description of the new batch sampling strategy in AMR. In Section~\ref{sec:illustrative}, the performance of AMR-PBS is analyzed and compared with a Bayesian EGO algorithm, using a set of test benchmark problems. For studying the effectiveness of the proposed algorithm, a CFD-based problem which is computationally very expensive to simulate is considered. Section~\ref{sec:CFD} describes the automated CFD framework and its optimization formulation, and discusses the results. Finally, concluding remarks are given in Section~\ref{sec:conclusion}.

\section{Variable-Fidelity Optimization with Adaptive Model Refinement (AMR)} \label{sec:VFM-AMR}
In this section, we introduce \emph{Adaptive Model Refinement (AMR)}, which is a decision-making tool for the timing of surrogate model refinement and identifying the optimal batch size and location of infill points in surrogate-based optimization (SBO). Performing model refinement (i.e., increasing model fidelity) too early in the design process can be computationally expensive while wasting resources to explore undesirable regions of the design domain. On the other hand, updating the surrogate model too late might mislead the search process early on to sub-optimal regions of the design domain, that is, leading to scenarios where the global optimum is outside of the region spanned by the population of candidate solutions in later iterations. The implementation of the proposed AMR in a population-based algorithm involves the five major steps, which are illustrated within the AMR flowchart shown in Fig. \ref{fig:MethodologyG}. These five steps can be summarized as:
\begin{figure}[!htp]
\centering
  \includegraphics[width=0.90\textwidth]{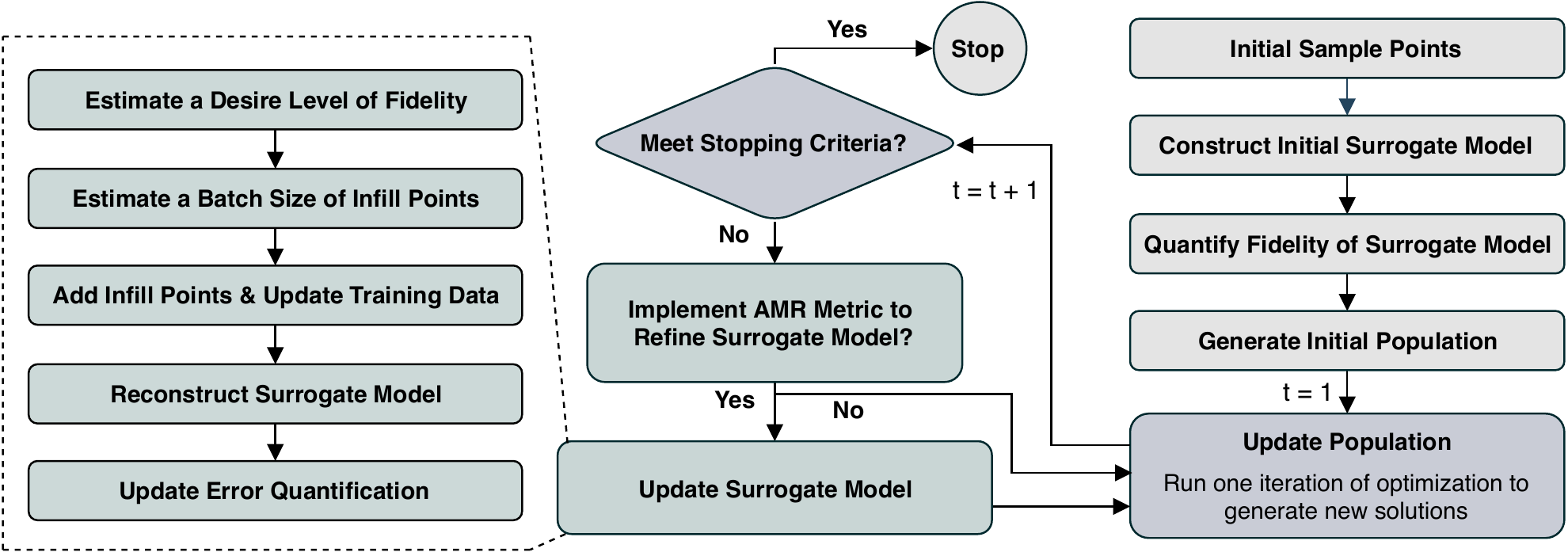}
  \caption{Adaptive model refinement (AMR) in surrogate-based optimization, with batch sampling}\label{fig:MethodologyG}
\end{figure}

\begin{description}
  \item[\textbf{Step 1}] This step occurs once. First, a set of initial sampling points are generated in the design space using a Latin hypercube (LH) sampling~\cite{McKay_79}. The best surrogate model is then constructed using the Concurrent Surrogate Model Selection (COSMOS) framework \cite{Chowdhury-COSMOS-SMO-2017,ghassemi2017optimal}, which automatically selects the best-suited model (from a pool comprising RBF, Kriging, SVR, and ANN models with different underlying kernel choices). The fidelity of the chosen surrogate model is given by PEMF \cite{mehmani2015predictive}, which provides a measure of the predicted error distribution. Finally, the initial population is generated at $t=1$, using the initial surrogate model. 
  \item[\textbf{Step 2}] At every iteration ($t$) of the heuristic optimization algorithm, the current surrogate model is used to update the function values of the population, and then set $t=t+1$. In this paper, particle swarm optimization is the chosen optimization algorithm.
  \item[\textbf{Step 3}] The stopping criterion is checked. The following two different methods can be used as the stopping criteria: (i) the difference between optimum values of five consecutive iterations is less than a threshold value, and (ii) the maximum allowed number of evaluations of function is reached. In this paper, the optimization algorithm stops when the relative change in the fitness function values in five consecutive iterations is less than a predefined function tolerance, $\delta_F$. If the termination criterion is satisfied, the current optimum  (the best global solution in case of PSO) is identified as the final optimum and the optimization process is terminated. Otherwise, go to \textbf{Step 4}.
  \item[\textbf{Step 4}] The AMR metric is evaluated in this step. If the AMR metric is satisfied, go to \textbf{Step 5}. Otherwise go to \textbf{Step 2}.
  \item[\textbf{Step 5}] In this step, a model refinement occurs, and the surrogate model is updated through a process that involves five sub-steps (see the expansion of ``Update Surrogate Model" step in fig. \ref{fig:MethodologyG}), of which the batch sampling is the important new addition. 
  Go to \textbf{Step 2}
\end{description}

In practice, the AMR metric (Step 4) is not necessary to be applied at every iteration; the user can specify that it be applied after every $\tau$ iteration. The first version of the AMR is presented in our previous work~\cite{mehmani2015adaptive, mehmani2015variable}, hence in this work we only focus on a new extension of the AMR. The details of Step 1 to Step 3 can be found in Appendix~\ref{app:amr}. In the next subsection, we briefly talk about Step 4 and then expand on the batch sampling aspects of Step 5. 

\subsection{The Adaptive Model Refinement (AMR) Metric}\label{ssec:metric}
In AMR, it is assumed that the uncertainty associated with surrogate model can be evaluated in the form of an error distribution, $\mathbb{P}_i$. Under this assumption, the fitness function values evaluated using the current surrogate model (in the SBO process) can be related to the corresponding high-fidelity model estimation or experiment-based system evaluations as
\begin{equation}\label{Eq:ErrorMode}
Y(x)=\hat{f}(x)+\varepsilon
\end{equation}
Here, $\hat{f}(.)$ and $\varepsilon$, respectively, represent the response of the current surrogate model and the stochastic error associated with it, and $Y$ is the corresponding high-fidelity model (or experiment-based system) response.

The AMR metric is then defined based on  "\emph{whether the uncertainty associated with a surrogate model response is higher than the observed improvement in the relative fitness function of the population}". The proposed model refinement metric is designed to use the stochastic global measures of surrogate model error and the distribution of solution improvement. This metric is formulated as a statistical test for surrogate model (SM) as follows
\begin{equation}
\begin{split}
& H_{0}:~\mathbb{Q}_{\mathbb{P}_{_{SM}}}(p_{cr}) \geq \mathbb{Q}_{\Theta}(1-p_{cr})
\end{split}
 \label{Eq:HTest}
\end{equation}
\noindent where $p_{cr}$ is the critical probability ($0<p_{cr}<1$) and $\mathbb{Q}$ represents a quantile function of a distribution.

\begin{figure*}[!htp]
\centering
\subfigure[Rejection of the test; \textbf{Don't Refine a model}]
{
\includegraphics[height=0.25\textwidth]{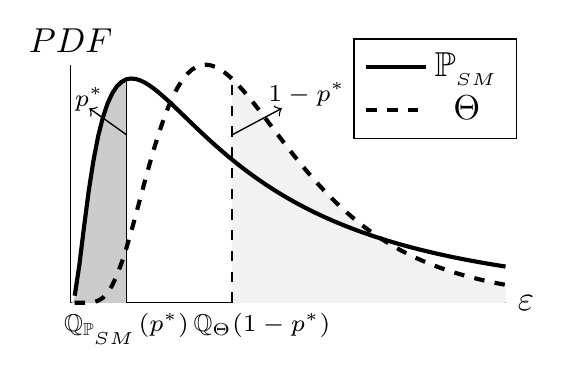}
\label{fig:amr_hypothesis_reject}
}
\subfigure[Acceptance of the test; \textbf{Refine a model}]
{
\includegraphics[height=0.25\textwidth]{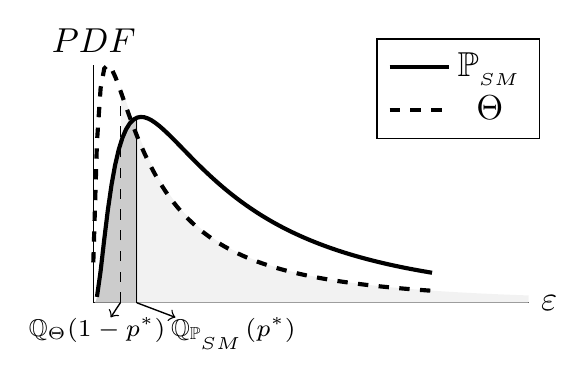}
\label{fig:amr_hypothesis_accept}
}
\caption{Illustration of the AMR hypothesis test (comparing the surrogate model error distribution ($\mathbb{P}_{_{SM}}$) and the distribution of fitness function improvement ($\Theta$)).}
\label{fig:amr_hypothesis}
\end{figure*}

For the sake of illustration, assume the surrogate model error distribution ($\mathbb{P}_{_{SM}}$) and the distribution of fitness function improvement ($\Theta$)) follow a log-normal distribution, and $p_{cr}=p^{*}$. In this case, the null hypothesis will be rejected, and the optimization process will use the current surrogate model \textbf{if}  $\mathbb{Q}_{\Theta}~>~\mathbb{Q}_{\mathbb{P}_{_{SM}}}$, as illustrated in Fig.~\ref{fig:amr_hypothesis_reject}. Conversely, \textbf{if}  $\mathbb{Q}_{\Theta}~<~\mathbb{Q}_{\mathbb{P}_{_{SM}}}$, the null hypothesis will be accepted, the surrogate model will be updated in the optimization process, as shown in  Fig.~\ref{fig:amr_hypothesis_accept}.

\subsection{Batch Sampling}
As mentioned above, the Adaptive Model Refinement (AMR) technique is a decision-making tool for the timing of surrogate model refinement, and for identifying not only the optimal batch size, but also the location of infill points needed to refine the surrogate model during a population based optimization process (e.g., genetic algorithm or PSO based process).

\subsubsection{Batch Size}
The strategy of tracking the variation of error with an increasing density of training points in the PEMF method~\cite{mehmani2015predictive} is then applied to determine the batch size for the samples to be added to the current sample points ($\Scale[0.95]{X^{CURR}}$). The inputs and outputs of PEMF (in the AMR method) can be expressed as:
\begin{equation}\label{Eq:PEMFfunction}
[ \mathbb{P}_{(\mu_{\varepsilon},\sigma_{\varepsilon})},~\Gamma^{Infill} ]=f_{_{PEMF}}(X,~\varepsilon^{*})
\end{equation}
where the vector $\rm X$ represents the sample data (input and output) used for training the surrogate model; and $\varepsilon^{*}$ is the desired fidelity in the model refinement process. In Eq.~\ref{Eq:PEMFfunction}, $\mathbb{P}_{(\mu_{\varepsilon},\sigma_{\varepsilon})}$, and $\Gamma^{Infill}$, respectively, represent the distribution of the error in the surrogate model, and the batch size for the infill points to achieve a desired level of fidelity in the AMR method. To this end, the desired batch size ($\Scale[0.95]{\Gamma^{Infill}}$) is estimated using the inverse of regression functions used to represent the variation of error with sample density (the Type 1 and Type 2 prediction functions used in PEMF), as shown below:
\begin{equation}\label{Eq:ModelRefinement3}
\mbox{$\Gamma^{Infill}$}=
\Bigg\{\left.
\begin{aligned}
\tiny
&\Scale[0.99]{\lceil \frac{\ln(\varepsilon_{mod}^{*})-\ln(a^{CURR})}{b^{CURR}} \rceil-N_s^{CURR}}&\mbox{\small Type~1}\\
&\Scale[0.99]{\lceil \exp (\frac{\ln(\varepsilon_{mod}^{*})-\ln(a^{CURR})}{b^{CURR}}) \rceil-N_s^{CURR}}&\mbox{\small Type~2} 
\end{aligned}
\right.
\end{equation}
where $\Scale[0.95]{a^{CURR}}$ and $\Scale[0.95]{b^{CURR}}$ are regression coefficients of the VESD function in the current surrogate model that are determined using the least square method, and $\Scale[0.95]{N_s^{CURR}}$ represents the size of the current sample points.  

\subsubsection{Batch Location}
The location of the new infill points ($\Scale[0.95]{\rm X^{Infill}| N(X^{Infill})=\Gamma^{Infill}}$) in the input space is determined based on a GP-based multi-points expected improvement criterion (q-EI). For a batch sampling, it is necessary to account for (and in our case mitigate) the interaction between the batch of future samples. Direct joint optimization of future candidate samples (in the batch)~\cite{azimi2010batch} is expensive to compute. Recently, Gonzalez et al.~\cite{gonzalez2016batch} reported a computationally tractable approximation to model the interactions, using a local penalization term. Based on this penalization factor, we adopt a criterion based on \textit{Expected Improvement with Local Penalization (EI-LP)}. The implementation of this approach is summarized in Algorithm~\ref{alg:LPqEI}.

In the given algorithm, we are calculating the penalty factor, $\gamma(x;x_i,\mathcal{L})$ (line 12 of Algorithm 1), which enables local exclusion zones based on the Lipschitz properties of the actual function ($f(x)$). This term tends to smoothly reduce the acquisition function in the neighborhood of the existing samples. To compute this penalty, we define a ball $\mathcal{B}_r$ with radius $\rho$ around each subsequently chosen sample point:
\begin{equation}
   \mathcal{B}_r(\mathbf{x},\mathbf{x}_i)
    = \{\mathbf{x}\in\mathcal{X}:\;\|\mathbf{x}_i - \mathbf{x}\| \leq \rho \}; ~
    \mathbf{x}_i\in X^{Infill}
\end{equation}
The local penalty associated with a point $\mathbf{x}$ is defined as the probability that $\mathbf{x}$ does not belong to the ball $\mathcal{B}_r$, as given by:
\begin{align}
\label{eq:localPenalizer}
    \gamma(\mathbf{x},\mathbf{x}_i) = 1 - P(\mathbf{x}\in \mathcal{B}_r)
\end{align}
By assuming that the distribution of the ball radius $\rho$ is Gaussian with mean $(M-\hat{f}(\mathbf{x}_i))/ L$ and variance $\sigma_r^2(\mathbf{x}_i)/L^2$, we can write the following expressions for the local penalty:
\begin{equation}
    \begin{aligned}
    \gamma(\mathbf{x},\mathbf{x}_i) &= 1 - P(\|\mathbf{x}_i - \mathbf{x}\|\leq \rho)\\
    &=P(\mathcal{N}(0,1)\leq\frac{L\|\mathbf{x}-\mathbf{x}_i\|-M+\hat{f}(\mathbf{x}_i)}{\sigma_r(\mathbf{x}_i)})\\
    &=\frac{1}{2}\text{erfc}\left(-\frac{L\|\mathbf{x}-\mathbf{x}_i\|-M+\hat{f}(\mathbf{x}_i)}{\sqrt{2\sigma_r^2(\mathbf{x}_i)}}\right)
    \end{aligned}
\end{equation}
Here, $M = \max_{\mathbf{x}} f(\mathbf{x})$ is the actual maximum value and $L$ is a valid Lipschitz constant ($\|f(\mathbf{x}_1)-f(\mathbf{x}_2)\|
\leq L\|\mathbf{x}_1-\mathbf{x}_2\|$).
\begin{algorithm}[!th]
\caption{Batch Sampling based on Local Penalizing}\label{alg:LPqEI}
\textbf{Input: $X^{CURR}, y^{CURR}, \hat{f}$}\\ \textbf{Output: $X^{UPDATED}$}
\begin{algorithmic}[1]
\State {$\mathcal{GP} \gets $ Fit a GP to $\mathcal{D}_n = \{X^{CURR}, y^{CURR}\}$}
\State {$u \gets (\hat{f}(x) - \min\{y^{CURR}\})/\sigma_n(x)$}
\State {$\alpha_{EI}(x; \mathcal{D}_n, \mathcal{GP}) \gets (u\Phi(u) + \phi(u))\sigma_n(x)$} \Comment{$\Phi(.)$: standard Gaussian distribution; $\phi(.)$: density functions.}
\State {$\mathcal{L} \gets \max_{x\in\mathcal{D}}\|\mu \Delta(x)\|$.} \Comment{$\mathcal{L}$ is a valid Lipschitz constant.}
\If {$\alpha_{EI}(x; \mathcal{D}_n, \mathcal{GP}) > 0$}
\State {$\hat{\alpha}_0(x) \gets \alpha_{EI}(x; \mathcal{D}_n, \mathcal{GP})$}
\Else
\State {$\hat{\alpha}_0(x) \gets \ln(1+e^{\alpha_{EI}(x; \mathcal{D}_n, \mathcal{GP})})$}
\EndIf
\State {$X^{UPDATED} \gets \emptyset$}
\For {$i = 1:N_s^{Infill}$}
\State {$x_i \gets \text{arg}\max_{x\in\mathcal{D}}\hat{\alpha}_{i-1}(x)$}
\State {$\hat{\alpha}_{i}(x) \gets \hat{\alpha}_0(x)\prod_{j=1}^i\gamma(x;x_i,\mathcal{L})$}
\State {$X^{UPDATED} \gets x_i \cup X^{UPDATED}$}
\EndFor
\end{algorithmic}
\end{algorithm}

The lower bound ($L^{\digamma}_j$) and the upper bound ($U^{\digamma}_j$) of the $\rm j^{th}$ dimension of the current sampling range ($\digamma$) are given by
\begin{gather}\label{Eq:ModelRefinement4}
\begin{split}
& L^{\digamma}_j=\min \{ x_j^{min},~x_j^{min} \}\\
& U^{\digamma}_j=\max \{ x_j^{max},~x_j^{max} \}
\end{split}
\end{gather}
where, $x_j^{min}$ and $x_j^{max}$ are, respectively, the lower and the upper bounds spanned by the current population of particles in the $\rm j^{th}$ dimension of the design space; $X_j^{min}$ and $X_j^{max}$ are the lower and upper bounds of the $\rm j^{th}$ dimension of the design space. 

\subsection{Optimization Algorithm: Particle Swarm Optimization}\label{ssec:pso}
The AMR can be integrated with most population-based optimization algorithms, such as evolutionary and swarm intelligence algorithms, where improvement can be tracked from one generation of the population to another. Particle swarm optimization (PSO) \cite{Kennedy-PSO-1995} is particularly well-suited to AMR. This is because identities of population members or search agents (known as particles) are preserved across generations, allowing ready-made tracking of solution quality improvement for each agent from one generation to another. 

In this paper, we use an advanced implementation of the PSO algorithm called mixed-discrete PSO (MDPSO), developed by \cite{Souma2013MDPSO}. Unlike the conventional PSO algorithm, MDPSO provides: (i) an explicit diversity preservation capability that mitigates the possibility of premature stagnation of particles, and (ii) an ability to deal with both discrete and continuous design variables. MDPSO has been used to solve a wide variety of highly non-convex (often multimodal) mixed-integer nonlinear programming problems in wind farm design \cite{chowdhury2013optimizing} and design of unmanned aerial vehicles \cite{Chowdhury-UAV-PF-JOA-2016}. Further description of the MDPSO algorithm can be found in the following paper~\cite{Souma2013MDPSO}.

\section{Benchmark Testing}\label{sec:illustrative}
\subsection{Testing Overview}
To evaluate the proposed AMR algorithm, two types of experiments are conducted using three distinct test benchmark functions, namely: 
1) the two-dimensional {\it Three-Hump Camel function}, 
2) the two-dimensional {\it Branin-Hoo function}, and 3) the six-dimensional {\it Hartmann6 function}. The detailed descriptions of these functions can be found in \cite{molga2005test}. In {\it{Experiment~1:}}, AMR is run to analyse its performance over the three test functions, and compare with that of Bayesian Efficient Global Optimization (BEGO) algorithm~\cite{tajbakhsh2013fully}. In {\it{Experiment~2:}} a batch size analysis is undertaken to explore the performance of AMR in the Hartmann6 test function, across multiple initial investment (batch size).
The results of the experiments are evaluated and compared in terms of \textit{relative absolute error} (RAE). The problem settings used (e.g., number of initial/total samples, number of runs per problem for robustness analysis, etc.) are given in Table \ref{tab:benchmark-settings}
\begin{figure}[!htp]
    \centering
    \subfigure[Three-Hump Camel]{\includegraphics[height=1.6in]{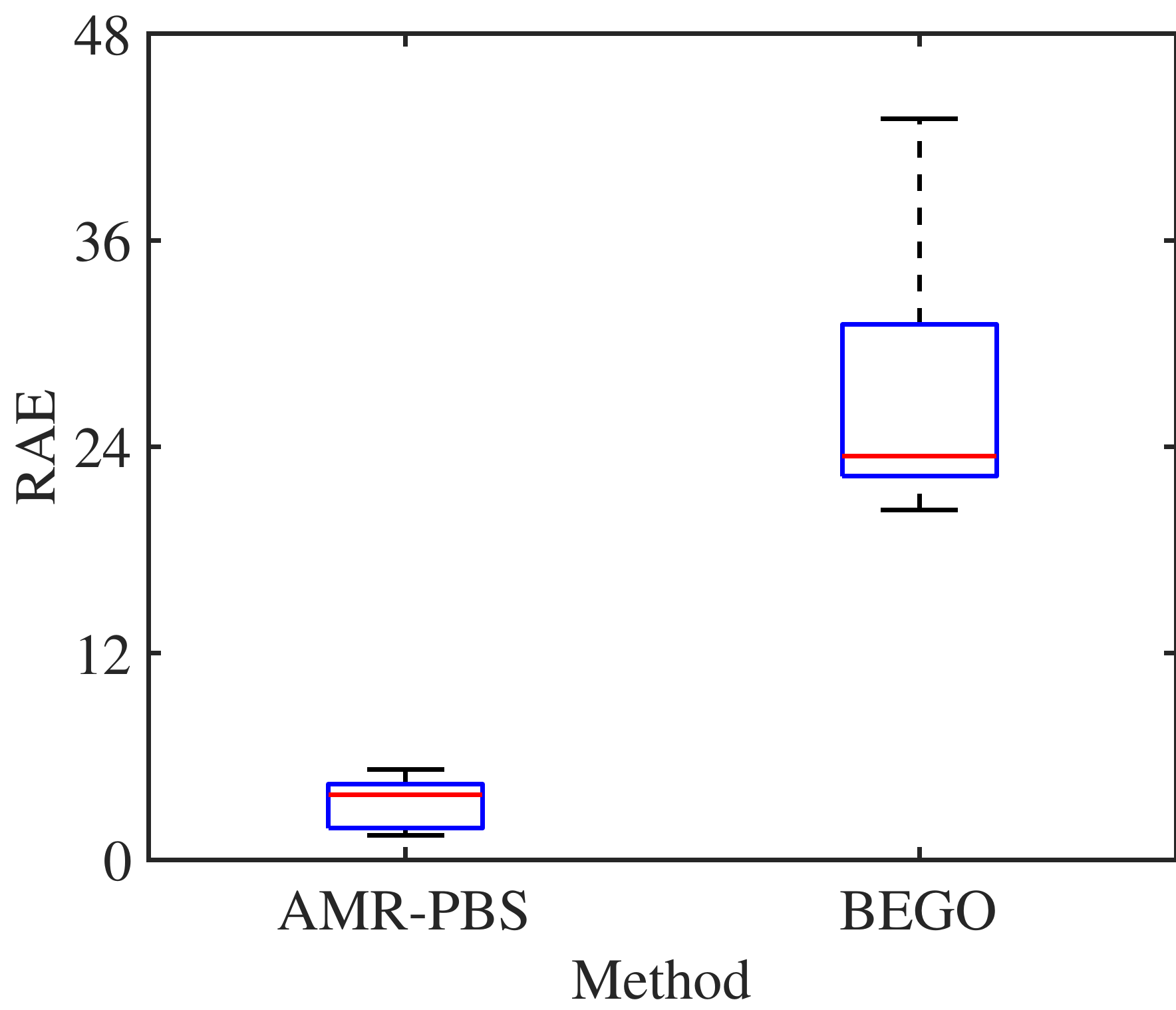}}%
    \hspace{0.2cm}
    \subfigure[Branin-Hoo]{\includegraphics[height=1.6in]{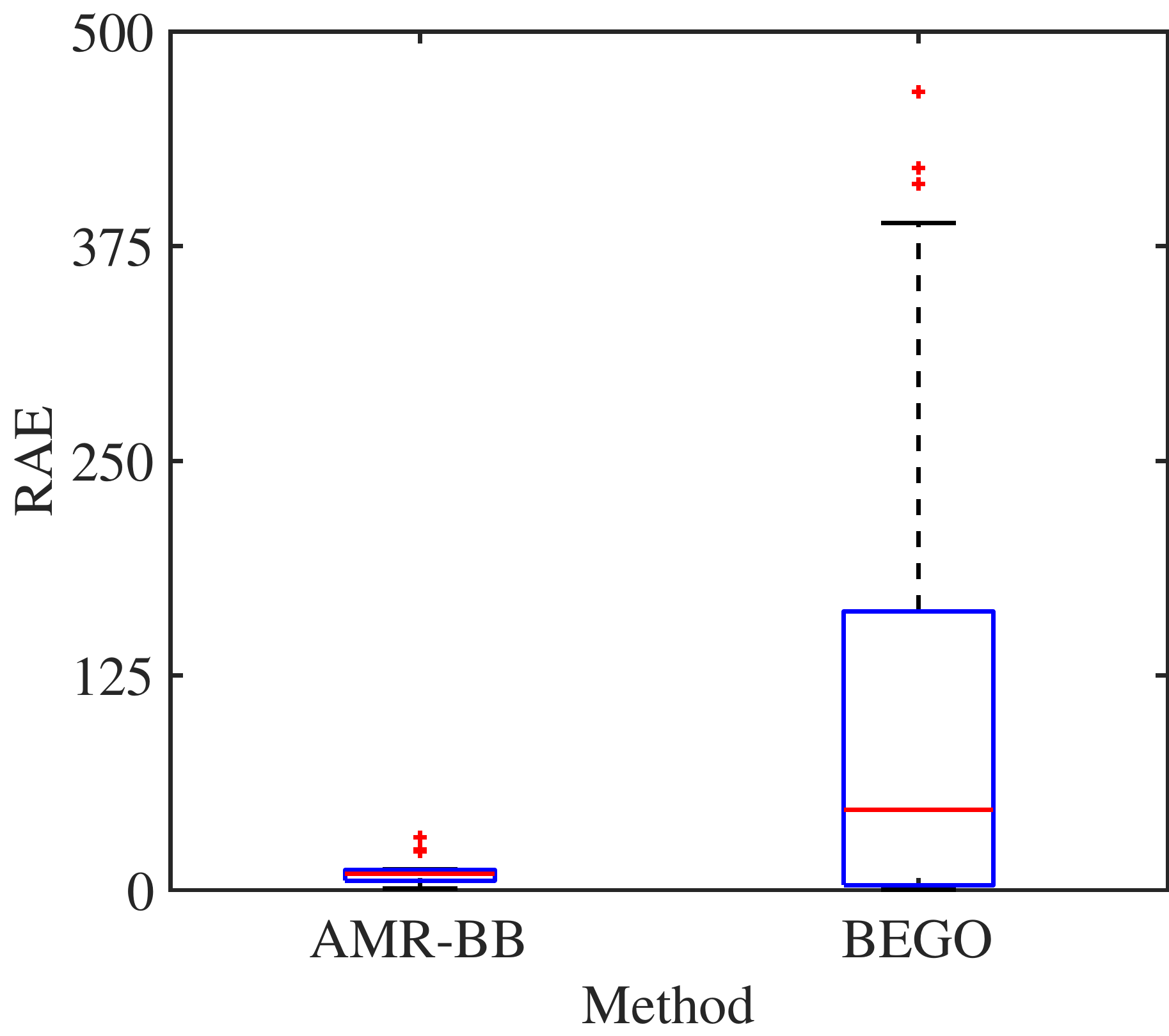}}%
    \hspace{0.2cm}
    \subfigure[Hartmann6]{\includegraphics[height=1.6in]{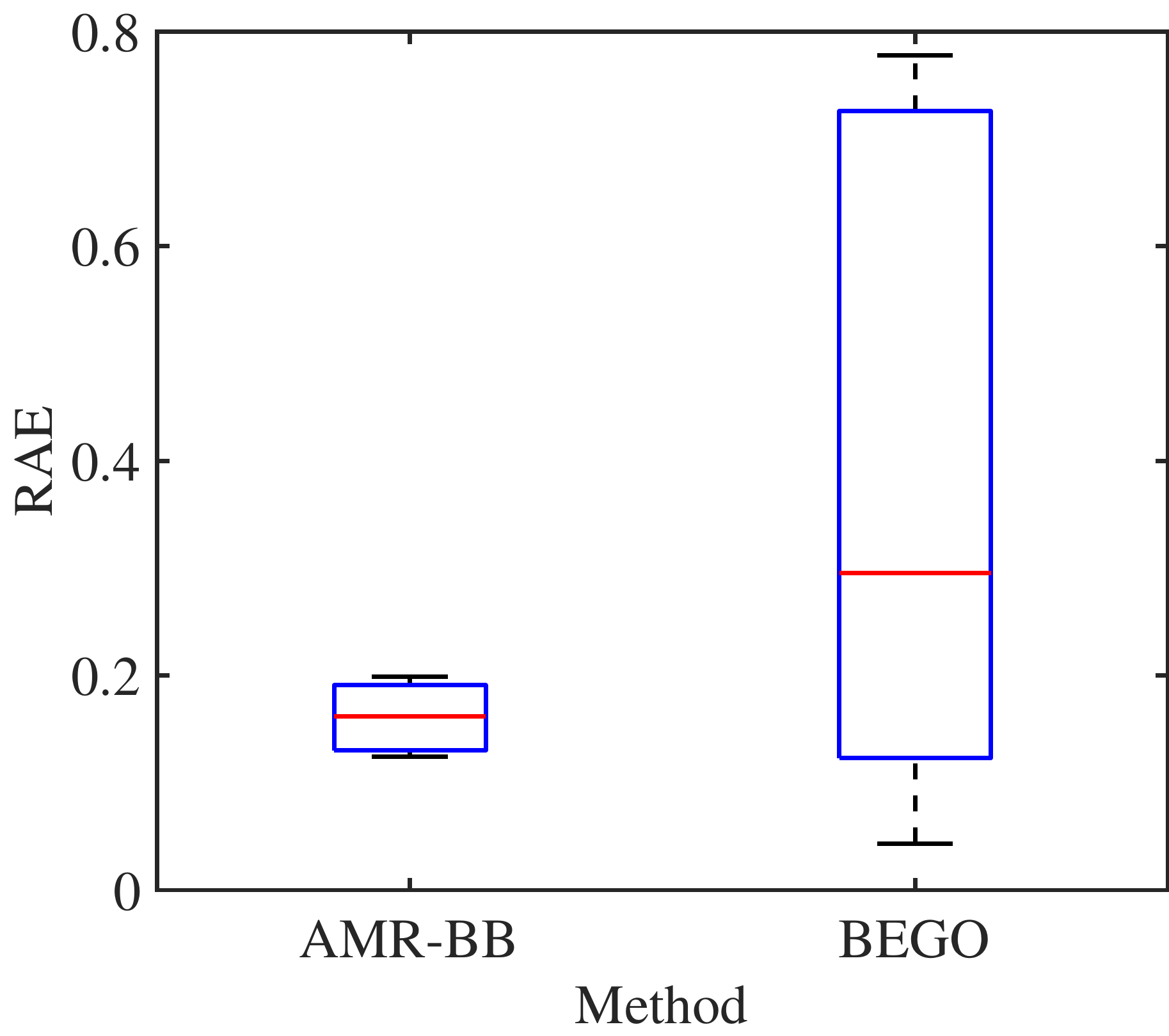}}
    \caption{Results on Test benchmark functions.}%
    \label{fig:AmrVsEgo}
\end{figure}

\begin{figure}
    \centering
    \subfigure[The convergence history and the surrogate model error; three refinements at iteration 16th and 24th with 2 and 8 infill points, respectively. The line with cross marker represents the surrogate model error. The line with circle shows the objective function (convergence history).]{
    \includegraphics[height=2in]{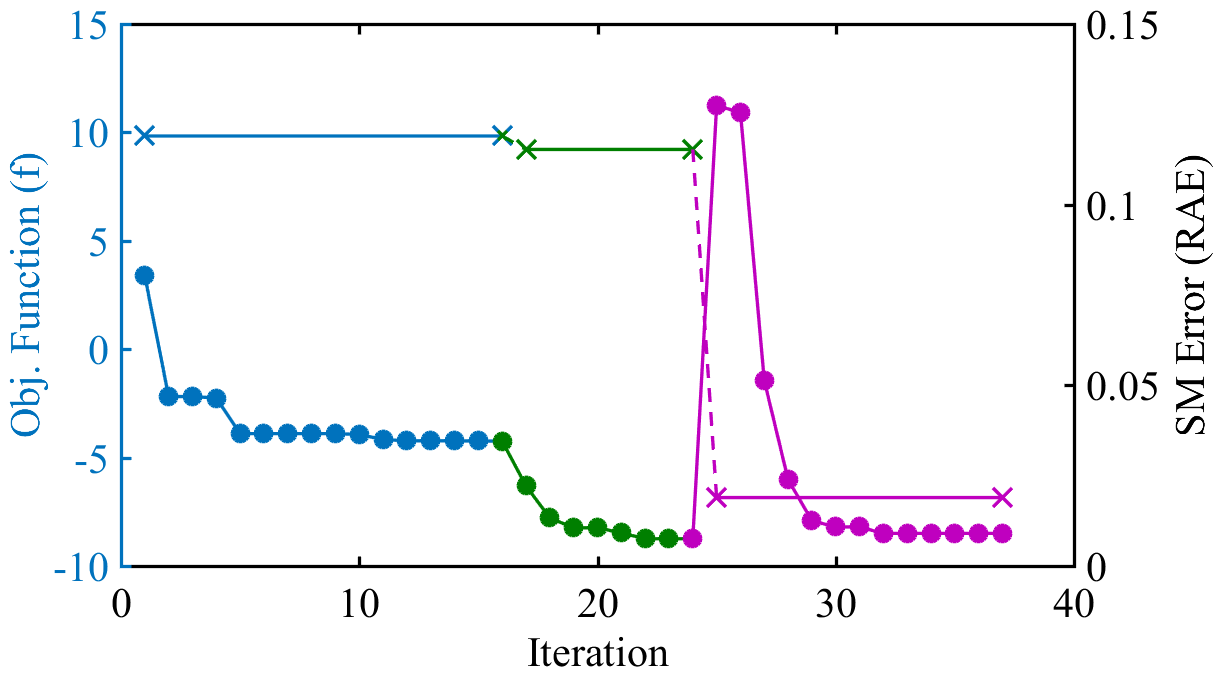}}\label{fig:convergenceHistory}%
    \hspace{0.2cm}
    \subfigure[The contour plot for the actual function and surrogate models for AMR. The solid line represent the actual function (ground truth). The dot and dash lines respectively represent the initial and final surrogate model in the AMR.]{\includegraphics[height=2in]{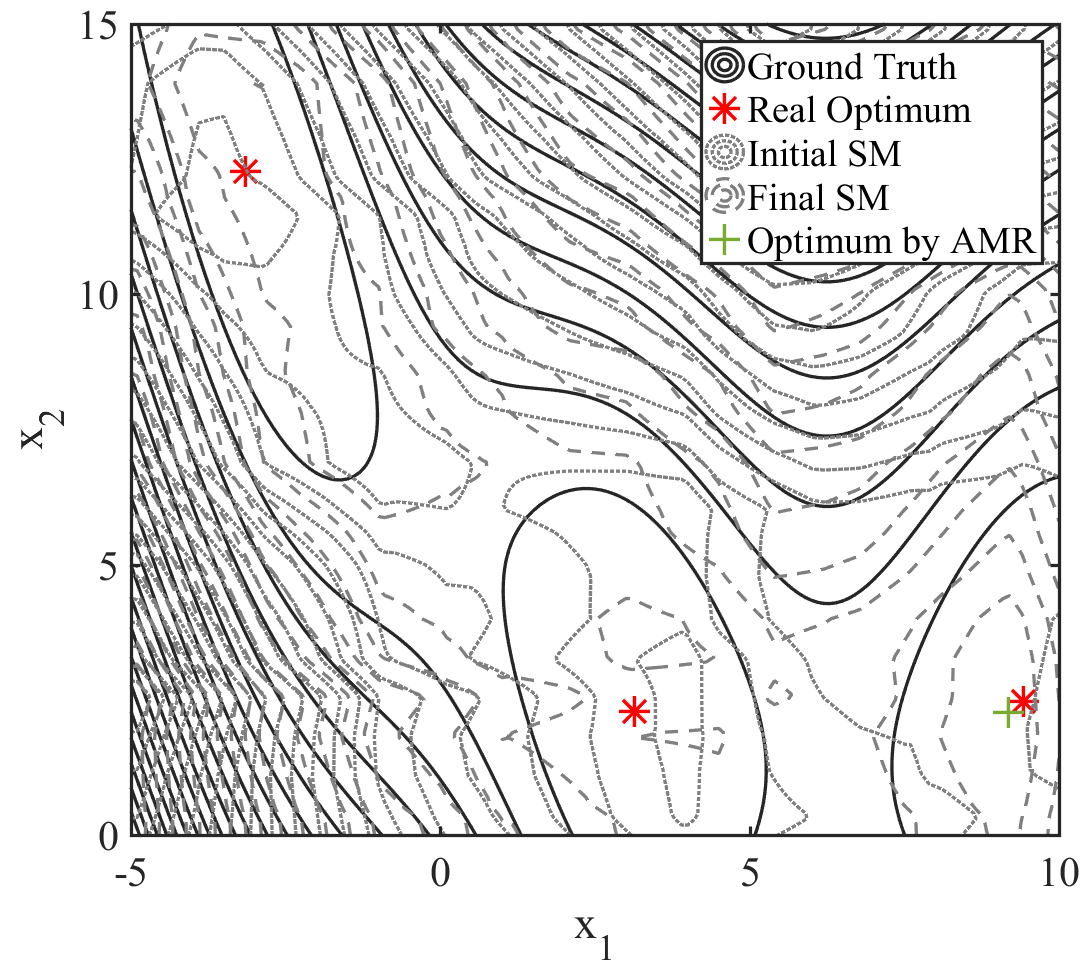}}\label{fig:responseSurface}
    \caption{The AMR results for the Branin-Hoo test function}
    \label{fig:selectedResult}
\end{figure}
\subsection{Experiment 1: A Comparative Analysis of AMR}
Figure~\ref{fig:AmrVsEgo} shows the performance of the AMR and the BEGO algorithms, in terms of RAE. As shown in these results, the AMR outperforms the BEGO algorithm in all the test functions, both in terms of the median value and variance across multiple runs (with randomized LHS samples). The median value and the variation of the results obtained by AMR is superior to BEGO, particularly, for Branin-Hoo and Three-Hump Camel functions. 

Figure~\ref{fig:selectedResult}~(a) shows the convergence history and the history of the surrogate model error over the optimization process. The contours for the actual function and the initial and final surrogate models in the AMR are illustrated in Fig.~\ref{fig:selectedResult}~(a). It can be seen from this plot that the AMR is able to reach to one of the optimums.

\begin{table}[ht!]
\begin{center}
\caption{Problem settings for the benchmark problems; $d$: problem dimension; $N_0$: initial sample size; $N_f$: final sample size.}
\begin{tabular}{ l c c c c c r}
\toprule[0.06em]
 \textbf{Benchmark Function} & $d$ & $N_0$  & $N_f$ & $N_\text{run}$ & \multicolumn{2}{c}{\textbf{RAE; Median (Max)}}\\
  & & & & & AMR-PBS & BEGO\\
\midrule[0.12em]
Three-Hump Camel & 2 & 20 & 30 & 20 & 4.39 (5.25) & 23.47 (43.07)\\
Branin-Hoo & 2 & 20 & 30 & 40 & 9.57 (30.63) & 46.68 (646.07) \\
Hartmann6 & 6 & 60 & 90 & 20 & 0.16 (0.19) & 0.30 (0.78) \\
\bottomrule[0.06em]
\end{tabular}
\label{tab:benchmark-settings}
\end{center}
\end{table}
\subsection{Experiment 2: A Parameter Analysis of Initial Sample Size}
\begin{figure}[!htp]
    \centering
    \includegraphics[width=0.45\textwidth]{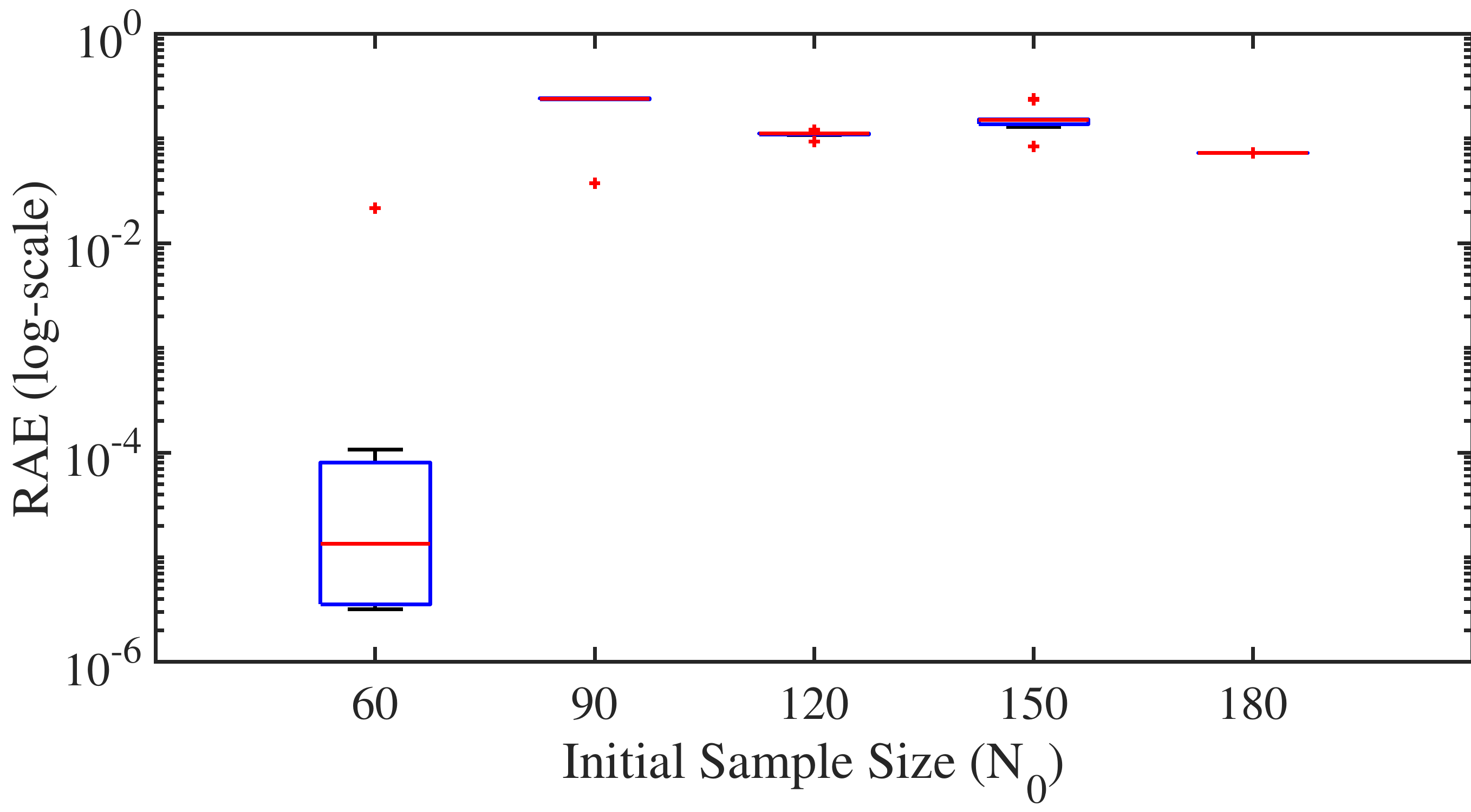}
\caption{Hartmann6 function: AMR performance variation with different initial/total sample ratios.}
    \label{fig:AmrBatchSize}
\end{figure}

In the second experiment, we study the impact of the initial investment/sample size ($N_0$) on the performance of the AMR technique. For this purpose, we use five different initial sample sizes varying from 60 to 180 with step size 30, and the maximum investment is set at 210 samples. For each case, we run the AMR approach 10 times, where the critical probability is set at 0.3 (i.e., $p_{cr}=0.3$) and $\tau$ is set at 2. Also, the population size ($N_\text{pop}$) and the maximum iteration ($\text{Iter}_\text{max}$) are set at 60 and 100, respectively. 
The obtained results are shown in Fig.~\ref{fig:AmrBatchSize}. It can be seen from this figure that small initial sample size ($N_0=60$ in this problem) improves the chance to lead the optimizer to the actual optimum.  For initial sample size larger than 30, the performance of the optimizer decreases (and for large values is almost similar, which can be because of the low complexity of the Hartmann6 problem). The reason for this observation is that when the initial investment ($N_0$) is small, the optimization is allowed to explore the design space and garner significant improvements before the refinement events get invoked.
%

\section{Application of AMR-PBS: Optimization of Riblet-based Flow Tailoring}\label{sec:cfdProblem}
\subsection{Riblet-based Flow Tailoring}
Lulekar et al. \cite{lulekar2018cfd} presented a new approach to represent and optimize surface riblets inspired by passive surface features observed in marine animals, for the purpose of reducing drag coefficients. While understanding of a broad range of riblet dimensions that could provide aerodynamic benefits existed, actual benefits might be realizable only in a very narrow range, which in the case of nature, is an outcome of millions of years of evolution. Hence, investigating (hypothesized) aerodynamic efficiency gains with minute surface riblets can become a challenging endeavor. The only ways to capture the complex flow mechanisms of interest induced by $O(10^{-1})$ millimeter scale riblets on a surface of a length scale, $O(10^2)$ mm, are to conduct high-resolution computational fluid dynamic (CFD) simulations or physical (wind tunnel) experiments. The prohibitive time/cost of both approaches limits the potential number of trials needed to identify the right desired riblet geometry. In this paper, we are using a set of open-source tools for running high-fidelity CFD simulation, where each run needs approximately 16 hours of computational time, executed on the UB CCR academic clusters using 2 compute nodes with Intel Xeon Processor E5-2660 (25M Cache, 2.20GHz), 16 cores per node, and 48GB RAM, indicating the expense of one single CFD simulation. One solution for overcoming this computational cost is minimizing the number of high-fidelity evaluations. For this purpose, we use AMR, where the surrogate model is trained over CFD or physical experiments, can address this challenge, and effectively lead us to the very small range of aerodynamically-promising riblet geometry. 

In order to exploit AMR, an automated modeling and design framework is needed that connects the various processes from design generation/CAD modeling to solution post processing. A first of its kind framework in the riblet design (for passive flow control) domain was developed by Lulekar et al. \cite{lulekar2018cfd}. Abrief descriptin of the riblet geometry and the automated framework are given next.

\subsection{The Geometry of Riblets}
\begin{figure}[!ht]
\begin{center}
\subfigure[]{\includegraphics[width=0.57\textwidth]{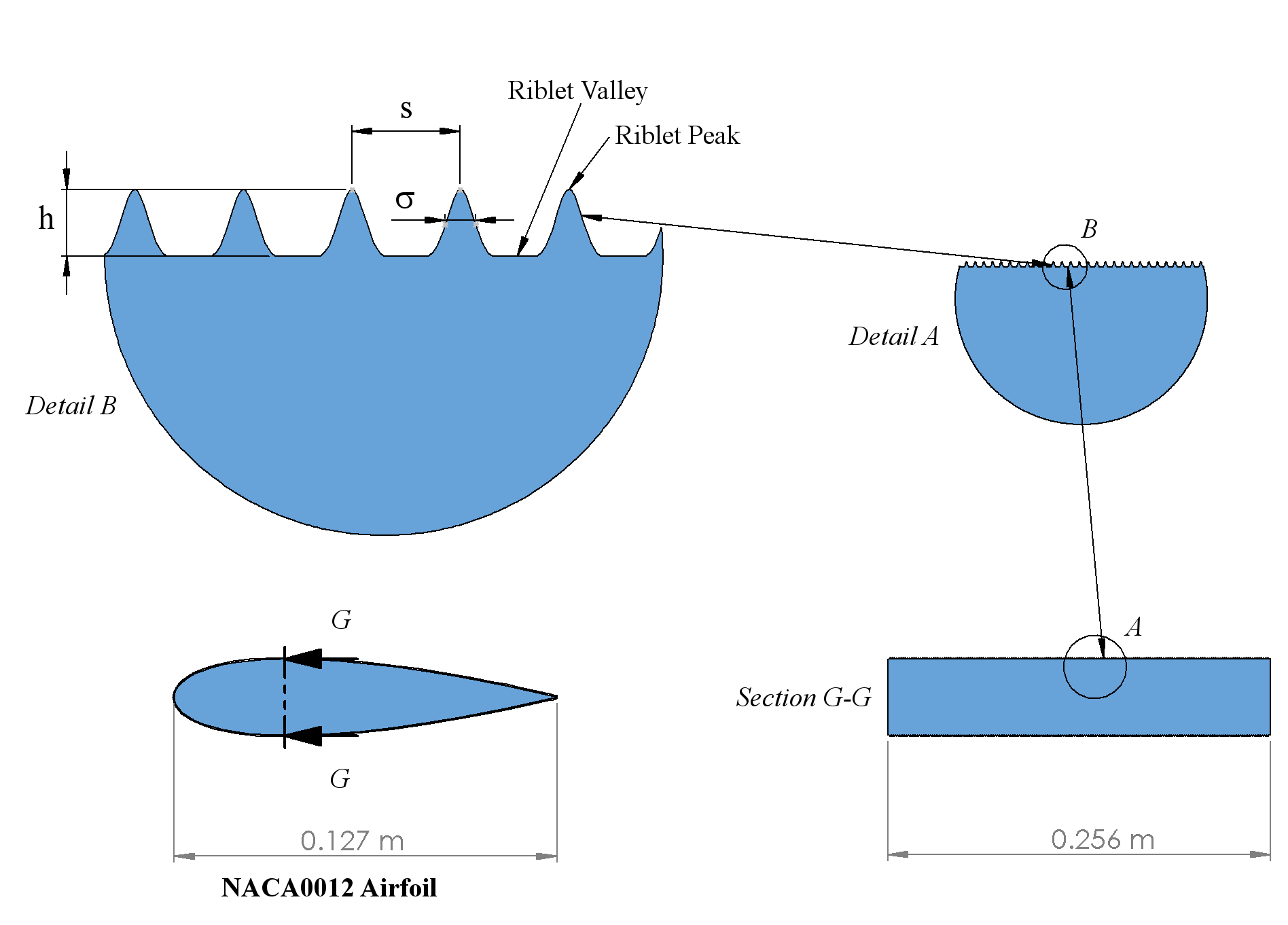}}%
\subfigure[]{\includegraphics[width=0.4\textwidth]{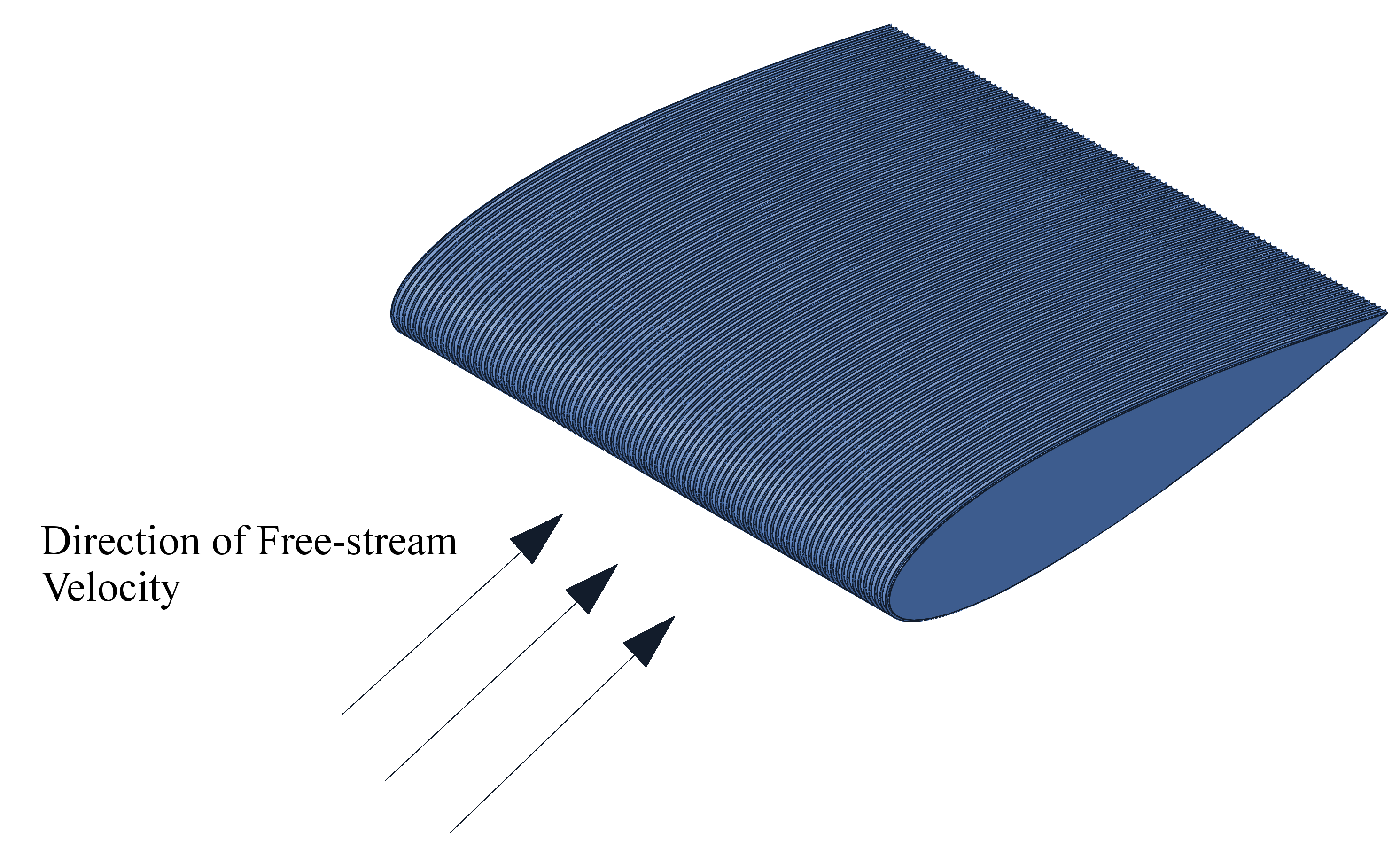}
  }
  \caption{(a) The shape and orientation of the Gaussian ridges; and (b) The ridgelines on a NACA0012 airfoil section, parallel to the free-stream velocity.}
  \label{fig:shapeANDorientation}
  \end{center}
\end{figure}

Different type of riblet shapes have been considered over decades to obtain drag reduction. 
We adopt our smooth riblet, which can be termed as patterns of Gaussian-shaped ridge lines~\cite{lulekar2018cfd}. Figure~\ref{fig:shapeANDorientation} shows the ridge-lines as a continuous extraction of a 2D Gaussian curve.  
The geometry of the ridge is parameterized in terms of peak height ($h$), spacing ($s$) and the distribution ($\sigma$) of the curve. Thus, any point $z$, on the ridgeline can be expressed as
\begin{equation} \label{eq:1}
z(y) = h e^{-{y^2}/{\sigma^2}}
\end{equation}
where $-3\sigma<y<3\sigma$. As illustrated in Fig.~\ref{fig:shapeANDorientation}, the pattern of ridgelines covers the entire top surface of the 3D airfoil section, and the cross-sectional shape/size of the ridge does not change in the streamwise direction in our current implementation. In most of the studies are conducted by changing the protrusion height ($h$) and the spacing ($s$) between the riblets to impede the streamwise vortices and their performance is measured based on a specific set of non-dimensional parameters, to account for the change in size of the flow structures like vortex diameter, where $h^+$ and $s^+$ is given as follows,
\begin{equation} \label{eq:2}
h^+ = \frac{hU_\tau}{\nu}\;\sqrt{\frac{C_f}{2}}
\end{equation}
\begin{equation} \label{eq:3}
s^+ = \frac{sU_\tau}{\nu}\;\sqrt{\frac{C_f}{2}}
\end{equation}
Further details of the physics and geometry of the riblets in the near wall flow, which affect the overall drag coefficient of the riblet surface can be found in our previous paper~\cite{lulekar2018cfd}. Here, instead, we focus on the modeling and optimization framework, which its central component is the presented AMR-based algorithm. A summary of this framework is described in the next.

\subsection{Automated CFD Framework}\label{sec:CFD}
\begin{figure}[th!]
\begin{center}
  \includegraphics[width=\textwidth]{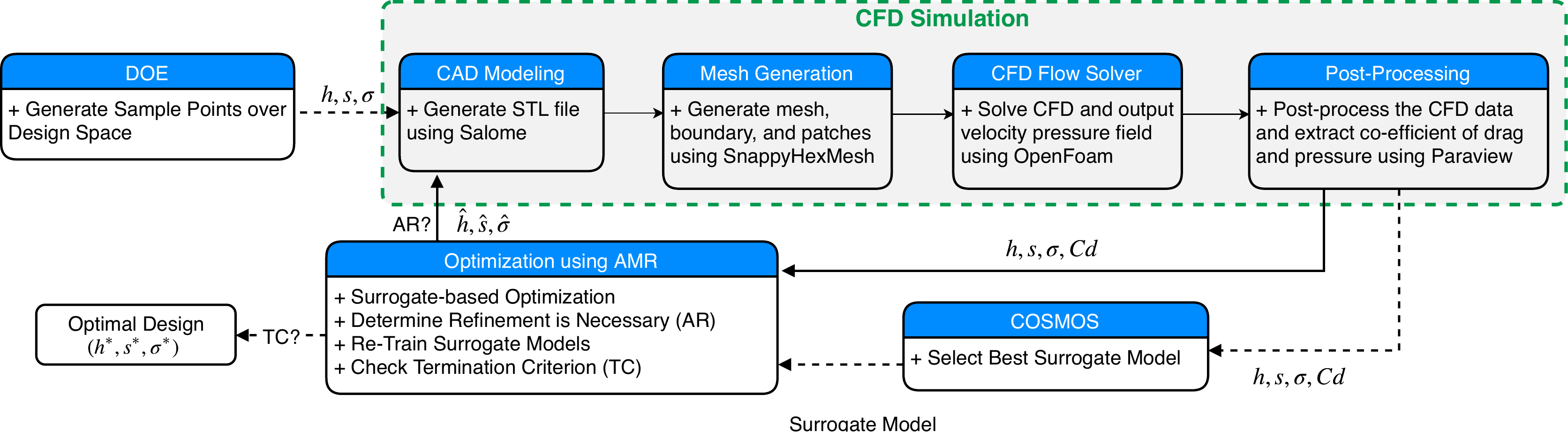}
  \caption{Overall optimization framework: integrating AMR with CFD simulation.}
  \label{fig:frameworkAMR}
  \end{center}
\end{figure}

In order to be able to study and solve the optimization problem to find the best drag reduction due to the riblets with the minimum human interaction, we developed an automated framework~\cite{lulekar2018cfd}, which integrates disparate computational tools in batch mode. The overall framework is illustrated in Fig.~\ref{fig:frameworkAMR} and it contains seven main blocks that can be divided into three categories: 1) Design-of-Experience (DoE) and metamodel training and selection (COSMOS) blocks; 2) CFD simulation blocks; and 3) multi-fidelity Optimization block. Each block is described in details in the next sections. 
Further descriptions of the CAD modeling, mesh generation, CFD flow solver, and post processing can be found in Appendix~\ref{app:cfdFramework}.

\subsection{Optimization Formulation and Design of Experiments}
The characteristics of the boundary layer are greatly affected by the type of the riblet used, spacing between two riblets and the protrusion into the boundary layer. Hence, identifying the suitable riblet shape boils down to three parameters, namely height, spacing and the distribution of the Gaussian curve. Given the objective of reducing drag ($D$) on the 3D airfoil section, the optimization problem can be stated as: 
\begin{equation} \label{eq:optimization}
\begin{aligned}
& \underset{\mathbf{x}=h,s,\sigma}{\text{minimize}}
& &  f = D(\mathbf{x}) \\
& \text{subject to}
& & g_1(\mathbf{x}) = 6\sigma - s \leq 0\\
& {}
& & g_2(\mathbf{x}) = s - 6h \leq 0\\
& {}
& & g_3(\mathbf{x}) = \sigma - 0.6h \leq 0
\end{aligned}
\end{equation}
where $h \in [0.2, 0.6]$, $s \in [0.72, 3.6]$ and $\sigma \in [0.12, 0.46]$, with the dimension of $h$ and $s$ are in millimeters. The first constraint, $g_1(\mathbf{x})$, is used to prevent adjacent Gaussian ridges from overlapping with each other. The second constraint, $g_2(\mathbf{x})$ in Eq.~\eqref{eq:optimization}, that restricts the inter-ridge spacing to 6 times the height of the ridge.
The third constraint, $g_3(\mathbf{x})$ in Eq~\ref{eq:optimization}, restricts the Gaussian curve to flatten out where the advantage of the riblets impeding the cross flow momentum is lost. 
Equations~\eqref{eq:2},~\eqref{eq:3},~\eqref{eq:optimization} can be used to set the bounds for $h$, $s$ and $\sigma$.

\subsection{Results}\label{sec:results}
\begin{figure}[th!]
\begin{center}
  \includegraphics[height=2in]{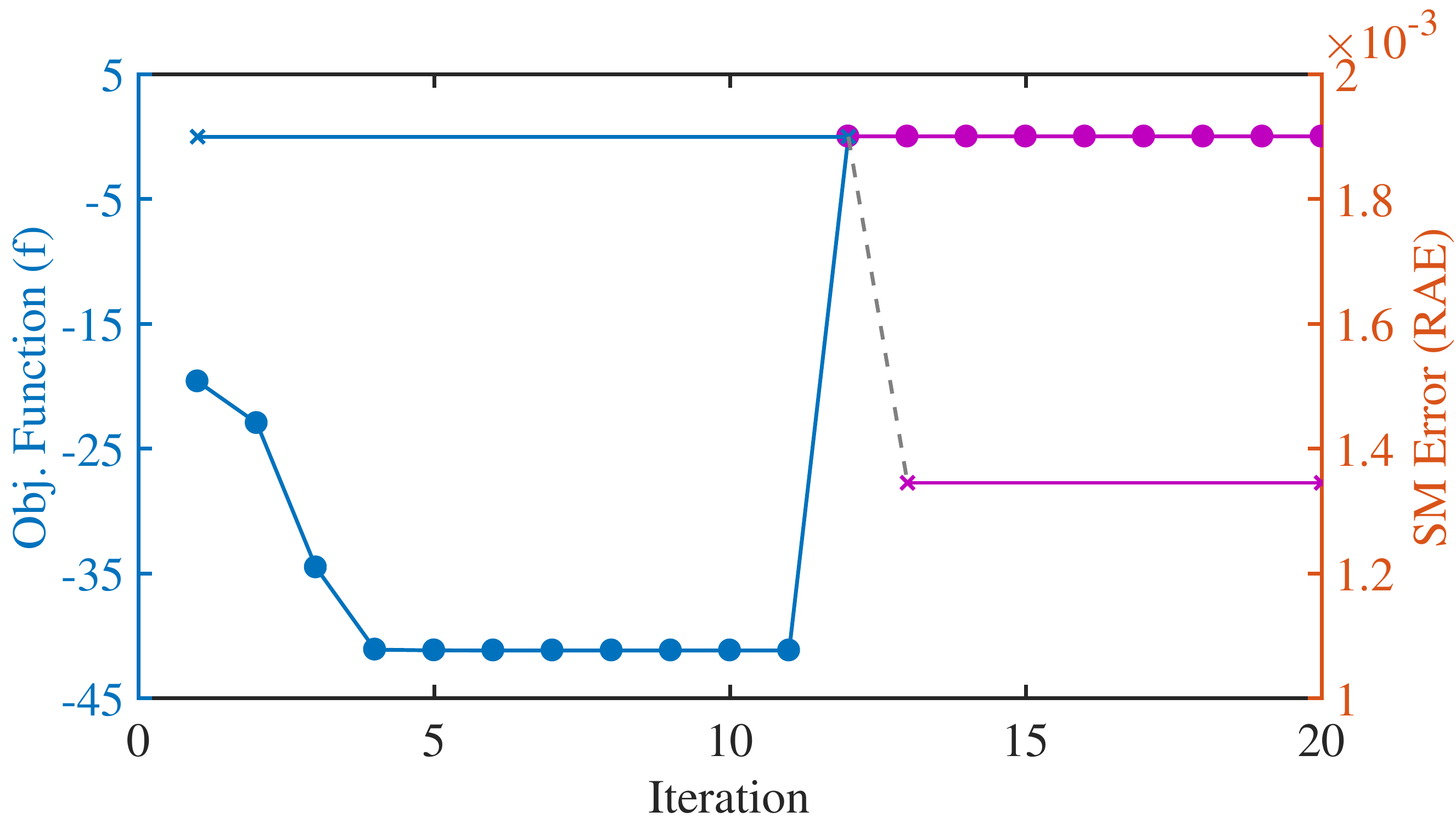}
  \caption{The convergence history and the surrogate model error for the CFD problem; As shown, it has only one refinement at iteration 12th with 29 infill points. The line with cross marker represents the surrogate model error. The line with circle shows the objective function (convergence history).}
  \label{fig:cfdAmrConvergence}
  \end{center}
\end{figure}

The simulation and optimization is performed for a case study with angle of attack $2^o$. This particular angle of attack is chosen, to observe the effect of riblet on drag reduction at higher pressure gradients. For this problem, the following settings are used: initial number of sample points (initial investment) used for training the surrogates is $N_0 = 30$, a population size of PSO is set at $N_\text{pop} = 30$, and the maximum iteration is set at $MaxIter = 100$. In addition, the critical probablity in AMR is set at $p_{cr} = 0.2$, and refinement events are allowed every 6th iteration during the optimization process. The convergence history of the AMR for this problem is shown in Fig.~\ref{fig:cfdAmrConvergence}. For this problem, the AMR requested only one refinement event at iteration 12th with 29 infill points. Table~\ref{tbl:AMR_Results} summarizes the results obtained for this problem. Kriging with Gaussian kernel is identified as the optimal surrogate model. The final results show that the optimal design variables (the chosen dimensions of the riblets) are $h^+=10.43$, $s^+=16.15$ and $\sigma=0.12$.
The coefficient of drag (Cd) at the optimum is 0.008046, which the surrogate model estimated it with 0.36\% accuracy. These dimensions provide a notable drag reduction (10.2\%) in compare with barefoil (airfoil without riblets).
\begin{table}[ht!]
\begin{center}
\caption{Results from Adaptive Model Refinement for $2^\circ$ Angle of Attack; $N_f$: final sample size; $N_\text{refin}$: number of refinements; $h^+$: Height; $s^+$: Spacing; $\sigma$: Distribution; RDC: Reduction in $C_D$ compared to baseline.}
\begin{tabular}{ l c c c c c c c c c c}
\toprule[0.06em]
 \textbf{SM} & $N_0$ & $N_\text{added}$ & $N_\text{refin}$ & $h^+$ & $s^+$ & $\sigma$ & \textbf{Cd-AMR} &\textbf{Cd-CFD} & \textbf{RAE [$\%$]} & \textbf{RDC [\%]}\\
\midrule[0.12em]
Kriging & 30 & 29 & 1 & 10.43 & 16.15 & 0.12 &  0.008017 & 0.008046 & 0.36 & 10.2\\
(Gaussian) &  & & & & & & & & &\\
\bottomrule[0.06em]
\end{tabular}
\label{tbl:AMR_Results}
\end{center}
\end{table}

\subsection{Flow Characteristics of the Optimum Design}
\begin{figure}[th!]
\centering
    \subfigure[The streamwise velocity.]{\includegraphics[width=0.48\textwidth]{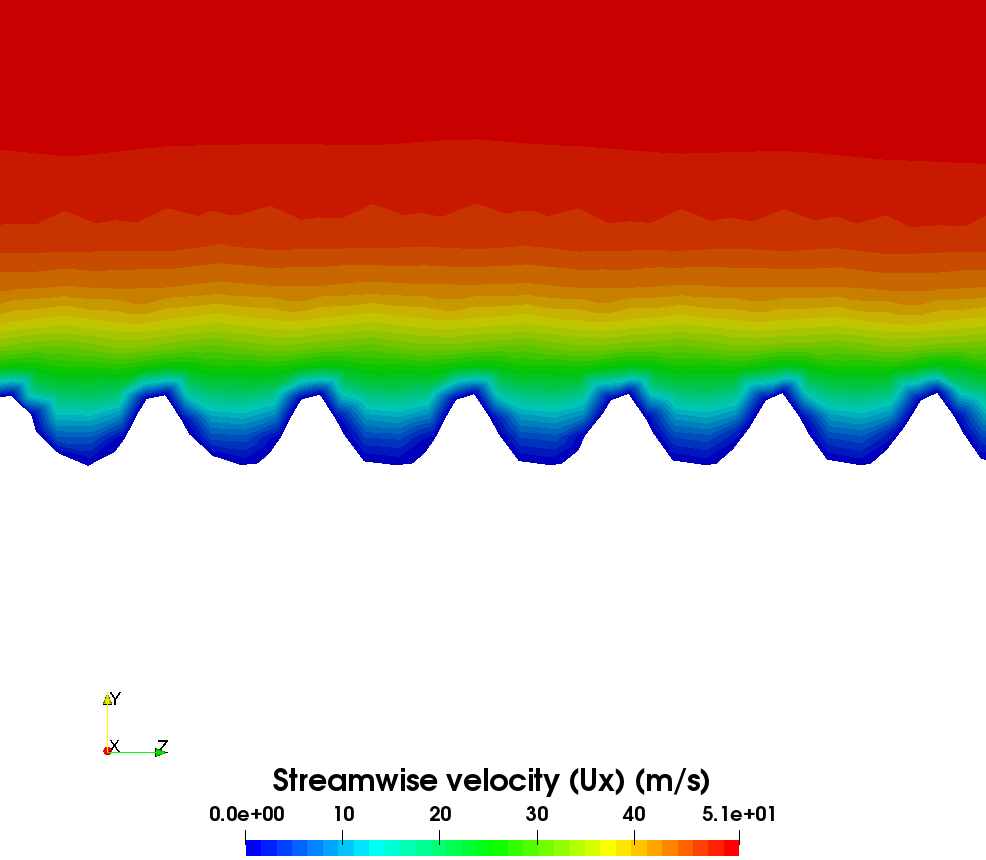}\label{fig:streamwise}}
  \hspace{0.2cm}
    \subfigure[The spanwise velocity.]{\includegraphics[width=0.48\textwidth]{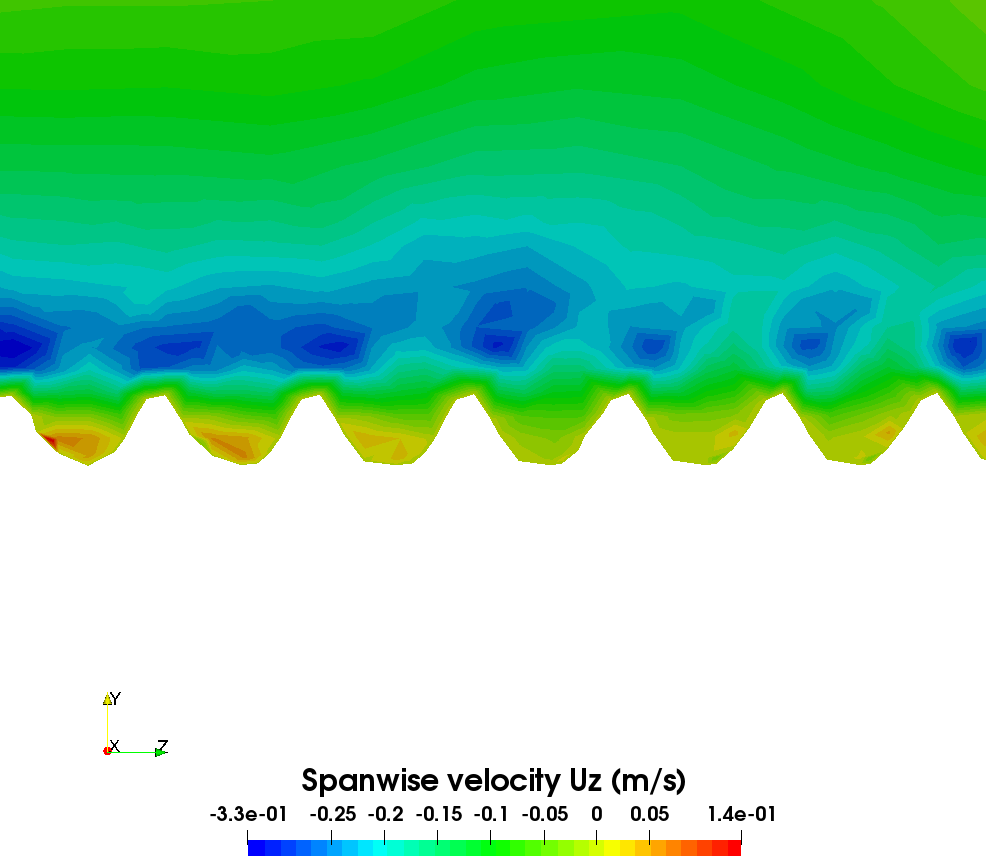}\label{fig:spanwise}}
    \caption{The contour plots at $x/c=0.3$ for the optimal design variables $h^+=10.43$, $s^+=16.15$ and $\sigma=0.12$.}
\end{figure}

Riblets are known to impede the cross-stream translation of the streamwise vortices. The flow very close to the wall explains the impending physics of the drag reduction. To understand the impending physics behind the reduction, we need to visualize the flow very close to the wall.
There is not a particular cross-stream velocity present but due to tip effect of the airfoil, a small amount of cross-flow stream is induced as visualized in Fig.~\ref{fig:spanwise}. Although a cross-flow is induced we can see a decrease in momentum transfer. It can be seen from Fig.~\ref{fig:streamwise} that the protrusion of the riblets from the surface into the boundary layer create zones of high-speed and low-speed fluid regions. The low-speed fluid region comes in contact with the riblet valleys which constitute a higher portion of the surface area whereas the high-speed fluid comes in contact only with the riblet peaks.

\begin{figure}[th!]
\centering
\subfigure[The Boundary Layer close to the wall compared with barefoil (zoomed)]{\includegraphics[width=0.45\textwidth]{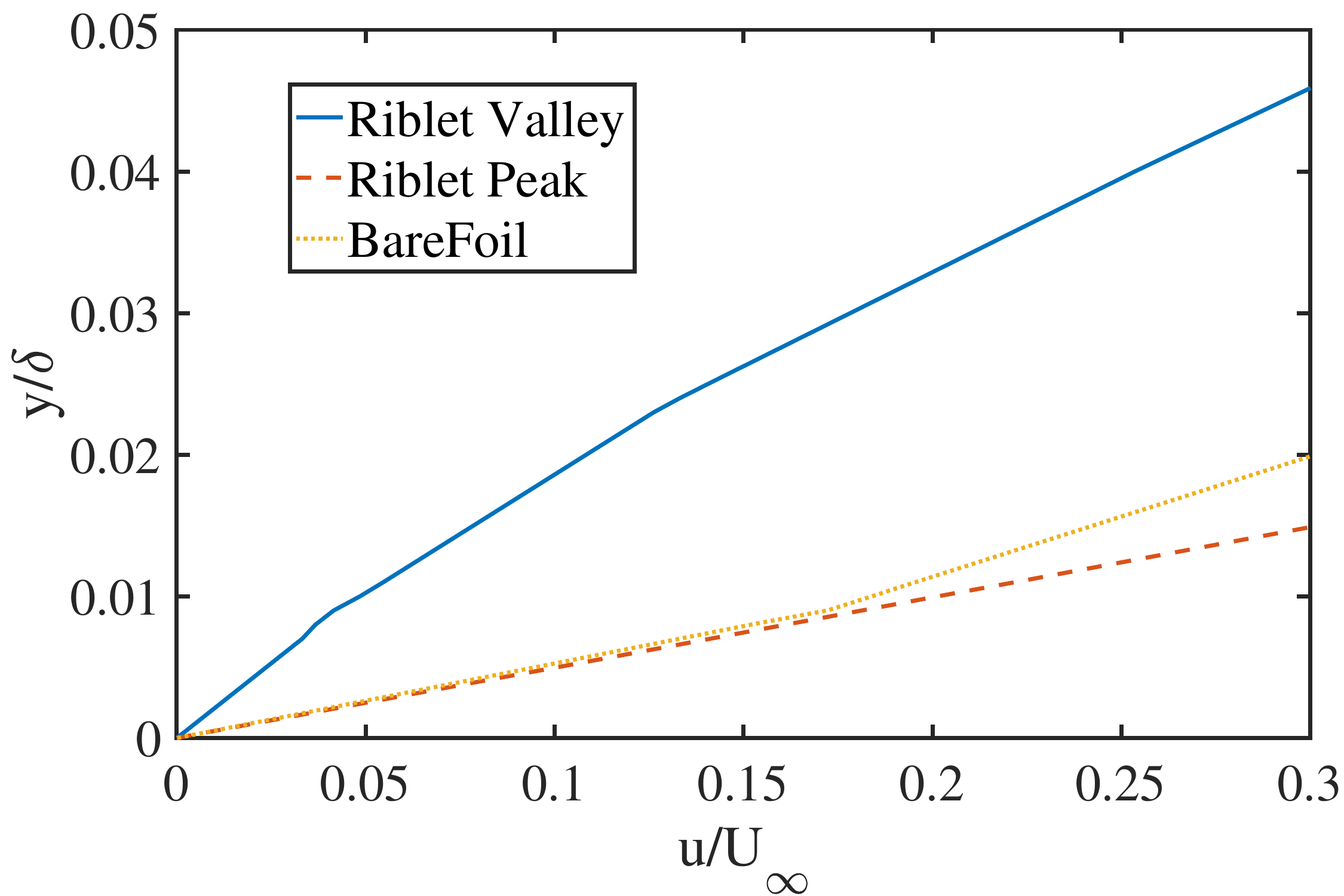}\label{optimum_30_boundary_layer}\label{fig:boundaryLayers}}%
\hspace{0.2cm}
\subfigure[The wall shear stress distribution over Gaussian riblets compare with the barefoil]{\includegraphics[width=0.45\textwidth]{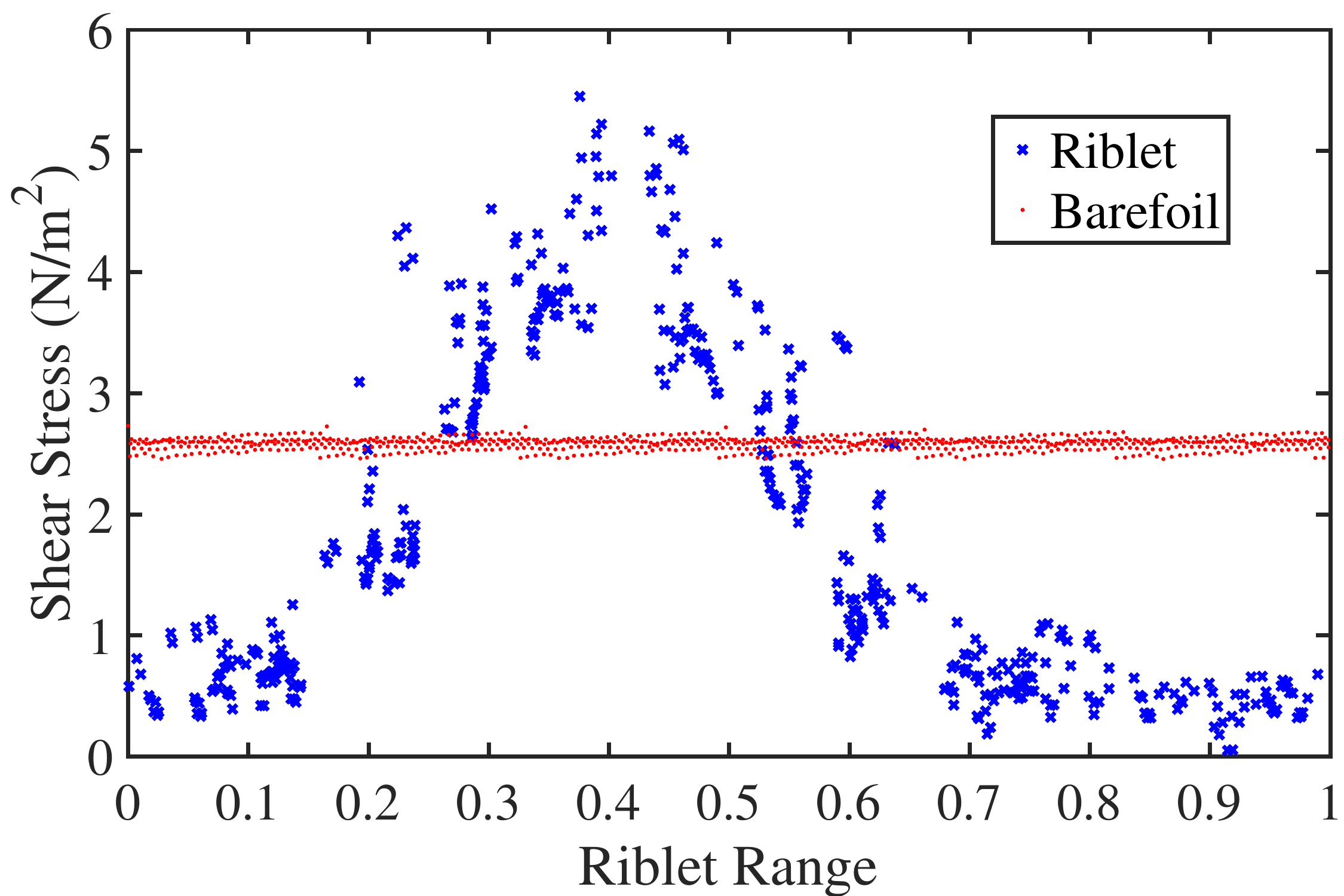}\label{fig:shearStress}}
\caption{A comparison between the optimal Gaussian riblets and the barefoil at $x/c = 0.3$}\label{fig:cfdResults}
\end{figure}

Figure~\ref{fig:boundaryLayers} shows the boundary layer in the ridge valley and at the ridge peak at $x/c=0.3$, compared with the boundary layer profile for the barefoil. This figure clearly shows the presence of different velocity gradients, attributed to low speed and high speed fluid regions. Velocity gradients are quantified by measuring the shear stress, where the wall shear stress is given by,
\begin{equation}
\tau_w = \mu \frac{dU}{dy}_{y=0} - \rho \overline{u'v'}
\end{equation}
Riblets demarcate the fluid into low speed and high speed regions, where the shear stress is simultaneously redistributed. In this study, the riblets are aligned with the free-stream velocity and hence, a shift in the wall shear stress distribution can be seen in the spanwise direction. Figure~\ref{fig:shearStress} shows the distribution of the wall shear stress when compared to the barefoil.

\section{Concluding Remarks}\label{sec:conclusion}
In this paper, we proposed a fundamental extension of a variable-fidelity optimization approach, so-called Adaptive Model Refinement (AMR), by incorporating a multi-point expectation improvement mechanism to decide the optimal location of batches of new samples. For evaluation, we used three benchmark functions and compared with an EGO-based approach (Bayesian EGO). The results show that the AMR-PBS performs 2-8 times better than the Bayesian EGO with significantly higher robustness.  A parameter analysis of initial sample size is performed, which shows that starting optimization with small sample size and utilizing resources once you need further accuracy can improve the optimization success. 

In addition, we applied the proposed variable-fidelity optimization (AMR-PBS) to design ribleted 3D airfoil surfaces for bio-inspired tailoring. While the results were not necessarily representative of converged results, since optimization were stopped due to time restrictions, the optimal design do show a promising drag reduction of 10.2\% in comparison with an airfoil without riblets. A large number of high and low speed regions are created due to the presence of the riblets balancing the adverse effects of increase in wetted surface area. Ideally, the riblets should be thin with sharp peak having minimal surface area to minimize the effect of high shear stress due to high-speed fluid contact. Although, a single high-resolution CFD simulation takes 7-16 hours, the given automated CFD framework, by utilizing the AMR technique and the parallel computing, was able to provision the optimal design in less than 40 hours. Future work entails investigation of how the interaction among samples in the batch impact the optimization process, and role played by the batch size in balancing potential interaction (and thus knowledge gain gap) and the need for reduction in model uncertainty for reliable continuation of the optimization process. Further application to an extended optimization of the passive flow surfaces to convergence will provide more comprehensive evidence regarding the effectiveness of the new AMR method. 

\bibliographystyle{IEEEtran}
\bibliography{cfdamr_reference}

\appendices 

\section{The Adaptive Model Refinement (AMR) Algorithm}\label{app:amr}
\subsection*{The Adaptive Model Refinement (AMR) Metric}\label{app:metric}
In this paper, it is assumed that the uncertainty associated with surrogate model can be evaluated in the form of an error distribution, $\mathbb{P}_i$. Under this assumption, the fitness function values evaluated using the current surrogate model (in the SBO process) can be related to the corresponding high-fidelity model estimation or experiment-based system evaluations as
\begin{equation}\label{Eq:ErrorMode2}
Y=\widehat{y}_{_{SM}}+\varepsilon
\end{equation}

In Equation~\ref{Eq:ErrorMode2}, $\widehat{y}_{SM}$ and $\varepsilon$, respectively, represent the response of the current surrogate model and the stochastic error associated with it, and $Y$ is the corresponding high-fidelity model (or experiment-based system) response. The relative improvement in the fitness function value ($\Delta f$) can be considered to follow an unknown distribution, $\Theta$, over the population of solutions. Here, $\Delta f$ in the $t^{th}$ iteration ($t\geq2$) can be expressed as
\begin{equation}
\mbox{$\Delta f^{~t}_{k}$}=
\Bigg\{\left.
\begin{aligned}
\small
|\frac{f^{~t}_k-f^{~t-\tau}_k}{f^{~t}_k}|  ~~&~ \mbox{if~$f^{~t}_{k}\neq0$} \quad\\
|f^{~t}_k-f^{~t-\tau}_k| ~~&~ \mbox{if~$f^{~t}_{k}=0$} \quad\\
 {\rm where}~k=1&,2,3,...,N_{pop}
\end{aligned}
\right.
\end{equation}
\noindent where, $\tau$ ($\in~\mathbb{Z}_{<t^*}$) is a user-defined parameter which regulates the occurrences of the ``surrogate model refinement" in the proposed SBO approach. Numerical experiments exploring different values for the $\tau$ indicate that the $3 \leq \tau \leq 5$ can be the suitable choice.

The AMR is then defined based on  "\emph{whether the uncertainty associated with a surrogate model response is higher than the observed improvement in the relative fitness function of the population}". The proposed model refinement metric is designed to use the stochastic global measures of surrogate model error and the distribution of solution improvement. Based on prior experience or practical design requirements, the designer is likely to be cognizant of what levels of global surrogate model error, $\eta$, is acceptable in an optimization process. Hence, $\eta$ can be perceived as a user preference. The \emph{critical probability, $p_{cr}$} for the current surrogate model (in the SBO process) with an error distribution $\mathbb{P}$ is then defined as the probability of the model error to be less than $\eta$. This definition can be expressed as

\begin{equation}\label{Eq:Pcr1}
p_{cr}=Pr(\varepsilon \leq \eta)=\int_{0}^{\eta}\mathbb{P}(\varepsilon ^{^‎\prime})~d\varepsilon ^{^‎\prime}
\end{equation}

The critical probability ($p_{cr}$) essentially indicates a critical bound in the error distribution $\mathbb{P}$ ($0 \leq \varepsilon \leq \eta$). If the predefined cut-off value ($\beta$) of the $\Theta$ distribution lies inside this region, the current surrogate model is considered to be no more reliable for use in the optimization process. 

The Adaptive Model Refinement (AMR) metric is formulated as a hypothesis for testing that is defined by a comparison between (a) the distribution of the relative fitness function improvement($\Theta$) over the entire population; and (b) the distribution of the error associated with the current surrogate model ($\mathbb{P}$) over the entire design space. This statistical test for one of the intermediate surrogate model (SM) in the SBO process can be stated as
\begin{equation}
\begin{split}
& H_{0}:~\mathbb{Q}_{\mathbb{P}_{_{SM}}}(p_{cr}) \geq \mathbb{Q}_{\Theta}(1-p_{cr}) 
\end{split}
 \label{Eq:HTest2}
\end{equation}
\noindent where $0<p_{cr}<1$ and $\mathbb{Q}$ represents a quantile function of a distribution; The $p$-quantile, for a given distribution function, $\Psi$, is defined as
\begin{equation}\label{Eq:quantile}
\mathbb{Q}_{\Psi}(p)=inf\{x \in \mathbb{R}:~p \leq \Psi_{(c.d.f.)}(x)\}
\end{equation}

In Eq.~\ref{Eq:HTest2}, the critical probability, $p_{cr}$, is an \emph{Indicator of Conservativeness (IoC)}. The IoC is based on user preferences, and regulates the trade-off between optimal solution reliability and computational cost in SBO. Generally, the higher the IoC (closer to 1), the higher the solution reliability and the greater the computational cost; under these conditions, model switching events will occur early on in the optimization process.

\begin{figure*}[!htp]
\centering
\subfigure[Rejection of the test; \textbf{Don't Refine a model}]
{
\includegraphics[height=0.25\textwidth]{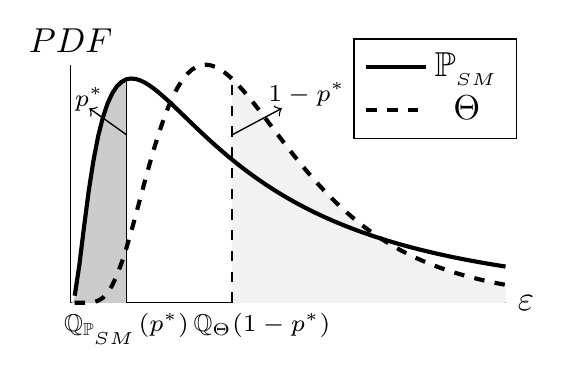}
\label{fig:amr_hypothesis_reject2}
}
\subfigure[Acceptance of the test; \textbf{Refine a model}]
{
\includegraphics[height=0.25\textwidth]{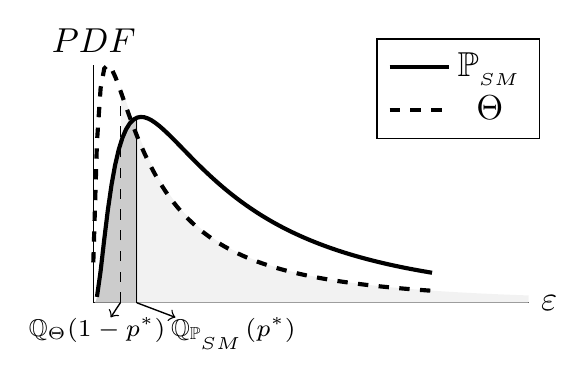}
\label{fig:amr_hypothesis_accept2}
}
\caption{Illustration of the AMR hypothesis test (comparing the surrogate model error distribution ($\mathbb{P}_{_{SM}}$) and the distribution of fitness function improvement ($\Theta$))}
\label{fig:amr_hypothesis2}
\end{figure*}

For the sake of illustration, assume $\Theta$ and $\mathbb{P}_{_{SM}}$ follow a log-normal distribution, and $p_{cr}=p^{*}$. In this case, the null hypothesis will be rejected, and the optimization process will use the current surrogate model \textbf{if}  $\mathbb{Q}_{\Theta}~>~\mathbb{Q}_{\mathbb{P}_{_{SM}}}$, as illustrated in Fig.~\ref{fig:amr_hypothesis_reject2}. Conversely, \textbf{if}  $\mathbb{Q}_{\Theta}~<~\mathbb{Q}_{\mathbb{P}_{_{SM}}}$, the null hypothesis will be accepted, the surrogate model will be updated in the optimization process, as shown in  Fig.~\ref{fig:amr_hypothesis_accept2}.
In this paper, kernel density estimation (KDE) is adopted to model the distribution of the relative improvement in the fitness function over consecutive ${\rm \Delta t}$ iterations. The uncertainty associated with surrogate models, and the batch size for the samples to be added in model refinement process of AMR are determined using an advanced surrogate error estimation method, called \emph{Predictive Estimation of Model Fidelity (PEMF)}~\cite{Ali_SMO_14}. The PEMF method is derived from the hypothesis that ``\emph{the accuracy of approximation models is related to the amount of data resources leveraged to train the model}''~\cite{Ali_SMO_14}. A brief description of the KDE and PEMF methods are provided in the following sections.

\subsection*{Model Refinement Based on PEMF}\label{ssec:ModelRefinemn}
In the previous section, the formulations and concepts of the AMR metric are defined. Based on the AMR metric, the model refinement will be performed at the $t^{*}$-th iteration of SBO under the condition
\begin{equation}\label{Eq:ModelRefinement2}
\mathbb{Q}_{\mathbb{P}_{_{\rm SM^{CURR}}}} \geq \mathbb{Q}^{t=t^*} _{\Theta}
\end{equation}
where $\rm SM^{_{CURR}}$ represents the current surrogate model in the optimization process. Model refinement is performed to efficiently improve the fidelity of the current surrogate model to achieve the desired fidelity for the upcoming iterations of SBO. In this paper, the desired fidelity, $\varepsilon_{mod}^{*}$, is determined using the history of the fitness function improvement in the optimization process, which is given by:
\begin{equation}\label{Eq:ModelRefinement1}
\varepsilon_{mod}^{*}=(1-\frac{\mathbb{Q}^{t=t^*}_{\Theta}-\mathbb{Q}^{t=t^*-\tau}_{\Theta}}{\mathbb{Q}^{t=t^*}_{\Theta}})~\times~\varepsilon_{mod}^{CURR}
\end{equation}

\noindent where, $\varepsilon_{mod}^{CURR}$ is the predicted modal error value associated with the current surrogate model.

\section{Automated CFD Framework}\label{app:cfdFramework}
\subsection*{Geometry and Meshing}
A finite span of NACA$0012$ airfoil is considered with riblets placed on the top and bottom surface of the airfoil. The airfoil contains a constant chord with an aspect ratio of 2. The chord length is 0.127m (5 inches) and is visualized in Fig.~\ref{fig:shapeANDorientation} along with the riblets. The riblets are a continuous extraction of a 2D Gaussian curve, where the riblet valleys and the peaks are aligned in the flow direction~\ref{fig:shapeANDorientation}. Three-dimensional CAD model of symmetric NACA$0012$ airfoil is build using the $SALOME-8.3.0$ and imported as an unstructured triangulated surface (.stl) with a precision of $O(10^{-2}) mm$. The surface mesh is imported for volumetric grid generation carried out using $SnappyHexMesh$. To
promote a greater stability with reasonable mesh count, while capturing the boundary layer profile more accurately, hex cells are preferred. Ten inflation layers are used with a growth rate of 1.2 and the smallest cell height is determined from $y^+\in[5, 10]$.

\begin{figure}[!ht]
\centering
    \subfigure[The prism Layers around the NACA$0012$ airfoil.]{\includegraphics[width=0.4\textwidth]{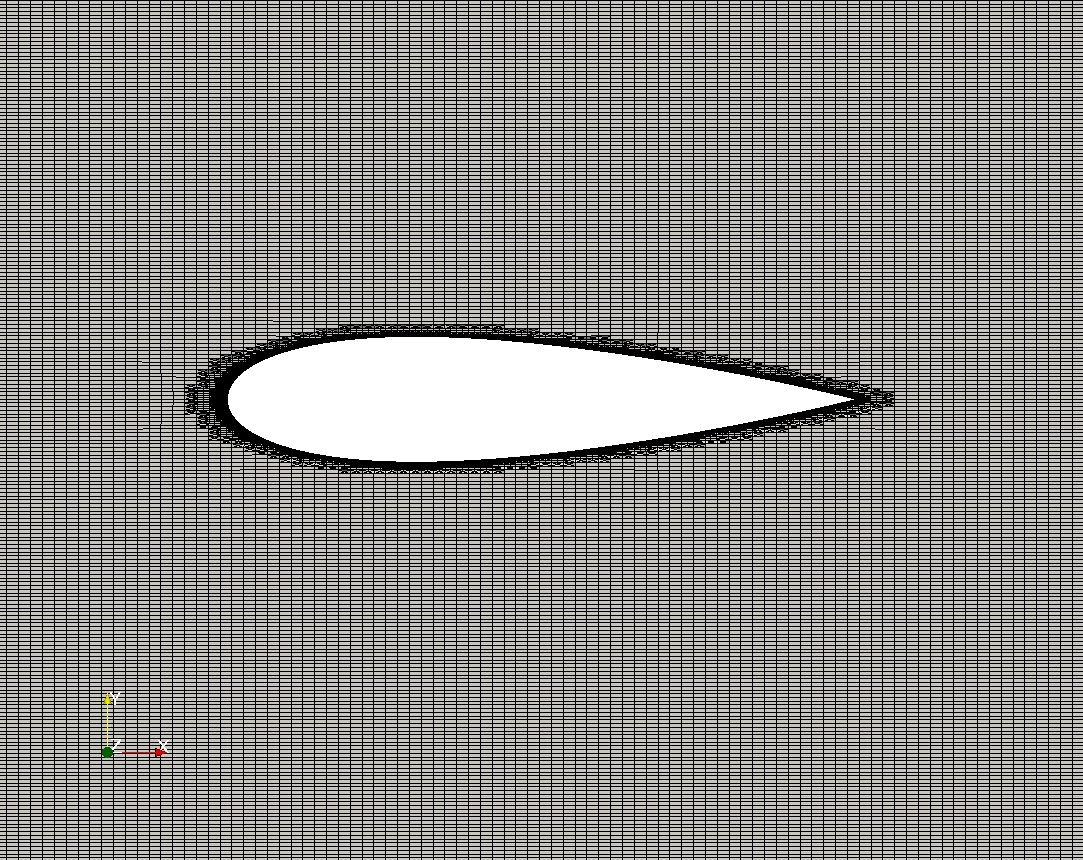}\label{fig:MeshAroundTheFoil}}
  \hspace{0.5cm}
    \subfigure[The mesh around the Gaussian riblets.]{\includegraphics[width=0.4\textwidth]{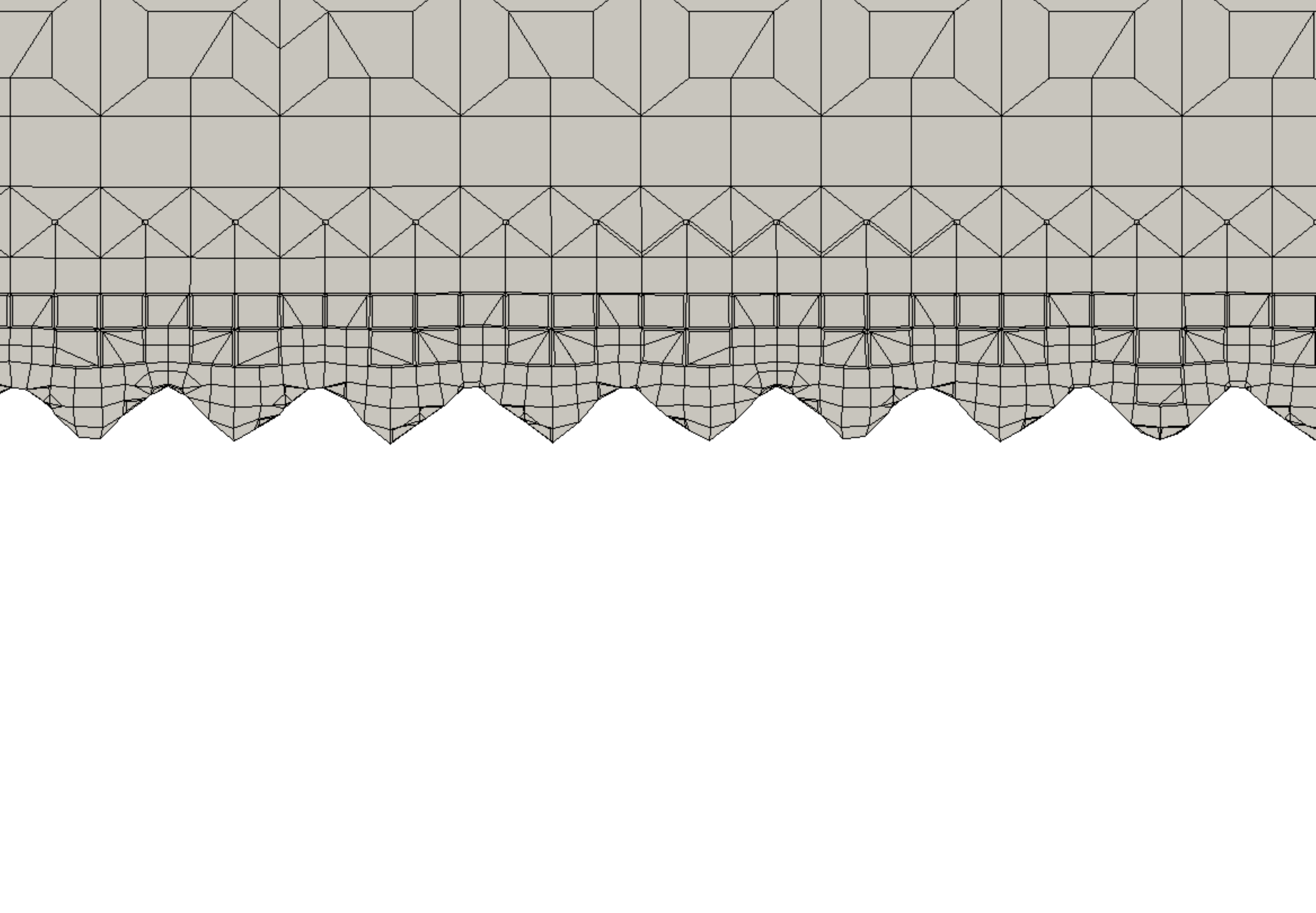}\label{fig:MeshAroundTheRiblet}}
   \caption{The prism layers and mesh around the airfoil and the Gaussian riblets.}
   \label{fig:MeshTheRiblet}
\end{figure}

To reduce the computational cost of the analysis, a symmetry plane is used and hence all the results are visualized on a half wing. The mesh generated in the domain approximately consists of 8 million cells. Fig.~\ref{fig:MeshAroundTheFoil} shows the prism layer around the NACA$0012$ airfoil and a sliced section of the mesh accurately capturing the Gaussian riblet curve can be seen in Fig.~\ref{fig:MeshAroundTheRiblet}.

\subsection*{Boundary Conditions}
The CFD simulations are carried out in a bounded box.
The domain has a velocity inlet ($free-stream$) with a constant velocity input. At the outlet, $pressure outlet$ for the pressure is provided and $zeroGradient$ for the velocities. At the top, bottom and side wall we have $slip$ condition for the velocities and $zeroGradient$ for the pressure. A $symmetry plane$ boundary condition is used to reduce the computational time, and $noSlip$ boundary conditions is applied on the airfoil. The flow being transitional region, a \%3.2 of turbulent intensity is given at the inlet. The turbulent kinetic energy (TKE) and the omega ($\omega$) at the inlet are calculated using the following equations:
\begin{equation}
I = 0.16Re^{-1/8}
\end{equation}
\begin{equation}
k = 1.5(UI)^2
\end{equation}
\begin{equation}
\omega = 0.9^{-1/4}\sqrt{k/l}
\end{equation}

\subsection*{Solver Settings}
An open source tool, $OpenFOAM$, is used to perform the CFD simulations. A pressure-based solver is used with $k-omega$ SST~\cite{menter1994two} turbulence models for wall bounded wall flow. The incoming $Re=3.38x10^5$, is calculated with a chord as the characteristic length, which is in the transitional ranger from laminar to turbulent. The fluid considered here is air, and treated as incompressible, as the Mach number is less than $0.4$. The pressure and velocity are coupled using the SIMPLE scheme. For spatial discretization and to calculate gradient of velocity ($\nabla U$), linear Upwind scheme is used. As stated above, the whole process is automated to evaluate multiple designs, \textit{OpenFoam} is used for CFD evaluations. 
Simulations are carried out with a termination criteria of $O(10^{-6})$. To shorten the computation time, \textit{MPI} library is used to parallelize the whole computation, and the simulations are performed using the distributed computing cluster (CCR) at the University of Buffalo. The primary quantities of interest obtained and derived from the CFD simulations are the shear stress data and co-efficient of drag.

\end{document}